\documentclass[12pt]{article} \usepackage{amsmath,theorem,
amssymb,amscd}

\input xypic

\newtheorem{conjecture}[equation]{Conjecture} \newtheorem{corollary}[equation]{Corollary} \newtheorem{lemma}[equation]{Lemma} \newtheorem{proposition}[equation]{Proposition} \newtheorem{theorem}[equation]{Theorem} \theoremstyle{plain}
\theorembodyfont{\rmfamily}  \newtheorem{definition}[equation]{Definition} 
 
 \newenvironment{pf}{\noindent \textbf{Proof.}}{\hfill{$\square$}\\} \numberwithin{equation}{section}

\begin{document}
\def\ch{\operatorname{ch}}
\def\td{\operatorname{td}}
\def\mod{\operatorname{mod}}
\def\gr{\operatorname{gr}}
\def\grmod{\operatorname{grmod}}
\def\A{\mathcal{A}}
\def\inf{\operatorname{inf}}
\def\discrep{\operatorname{discrep}}
\def\w{\omega}
\def\top{\operatorname{top}}
\def\O{\mathcal{O}}
\def\T{\Theta}
\def\m{\mathfrak{m}}
\def\F{\mathcal{F}}
\def\E{\mathcal{E}}
\def\faces{\mathcal{F}}
\def\P{\mathbb{P}}
\def\C{\mathbb{C}}
\def\N{\mathcal{N}}
\def\R{\mathbb{R}}
\def\Z{\mathbb{Z}}
\def\Q{\mathbb{Q}}
\def\ww{\wedge}
\def\even{\mbox{even}}
\def\odd{\mbox{odd}}
\def\a{\alpha}
\def\b{\beta}
\def\c{\gamma}
\def\d{\delta}
\def\D{\Delta}
\def\e{\varepsilon}
\def\l{\lambda}
\def\O{\mathcal{O}}
\def\t{\tau}
\def\s{\sigma}
\def\di{\partial}
\def\z{\zeta}
\def\la{\langle} 
\def\ra{\rangle}
\def\rtar{\rightarrow}
\def\faces{\operatorname{faces}}
\def\can{\operatorname{can}}
\def\rank{\operatorname{rank}}
\def\pic{\operatorname{pic}}
\def\symPic{\operatorname{symPic}}
\def\H{\operatorname{H}}
\def\coh{\operatorname{coh}}
\def\Bl{\operatorname{Bl}}
\def\g{\operatorname{g}}
\def\udim{\operatorname{udim}}
\def\order{\operatorname{ord}}
\def\height{\operatorname{ht}}
\def\det{\operatorname{det}}
\def\Gl{\operatorname{GL}}
\def\rad{\operatorname{rad}}
\def\biDiv{\operatorname{biDiv}}
\def\biPic{\operatorname{biPic}}
\def\locPic{\operatorname{locPic}}
\def\rtr{\operatorname{rtr}}
\def\Stan{\operatorname{Stan}}
\def\Trg{\operatorname{Trg}}
\def\B{\mathcal{B}}
\def\Cl{\operatorname{Cl}}
\def\H{\mathcal{H}}
\def\T{\mathcal{T}}
\def\M{\mathcal{M}}
\def\cok{\operatorname{cok}}
\def\Mod{\operatorname{Mod}}
\def\E{\mathcal{E}}
\def\ker{\operatorname{ker}}
\def\Hom{\operatorname{Hom}}
\def\shom{\mathcal{H}om}
\def\Set{\operatorname{Set}}
\def\Grp{\operatorname{Grp}}
\def\Ext{\operatorname{Ext}}
\def\Tor{\operatorname{Tor}}
\def\tors{\operatorname{tors}}
\def\End{\operatorname{End}}
\def\Pic{\operatorname{Pic}}
\def\Spec{\operatorname{Spec}}
\def\Proj{\operatorname{Proj}}
\def\Adj{\operatorname{Adj}}
\def\tr{\operatorname{tr}}
\def\Der{\operatorname{Der}}
\def\Sing{\operatorname{Sing}}
\def\OutDer{\operatorname{OutDer}}
\def\Obs{\operatorname{Obs}}
\def\AlgExt{\operatorname{AlgExt}}
\def\Exal{\operatorname{AlgExt}}
\def\Def{\operatorname{AlgExt}}
\def\kod{\operatorname{kod}}
\def\Open{\operatorname{Open}}
\def\supp{\operatorname{supp}}
\def\Br{\operatorname{Br}}
\def\EmbDef{\operatorname{EmbDef}}
\def\Aut{\operatorname{Aut}}
\def\Int{\operatorname{Int}}
\def\depth{\operatorname{depth}}
\def\lt{\operatorname{lt}}
\def\cal{\mathcal}
\def\L{\mathcal{L}}
\def\Sym{\operatorname{Sym}}
\def\cEnd{\cal{E}{nd}}
\def\Obs{\operatorname{Obs}}
\def\Out{\operatorname{Out}}
\def\ideal{\operatorname{I}}
\def\cir{\cdot}
\def\gkdim{\operatorname{gkdim}}
\def\GL{\operatorname{GL}}
\def\SL{\operatorname{SL}}
\def\Irr{\operatorname{Irr}}
\def\qcoh{\operatorname{qcoh}}
\def\id{\operatorname{id}}
\def\pd{\operatorname{pd}}
\def\G{\Bbb{G}}

\title{Conic Bundles and Clifford Algebras}
\author{Daniel Chan
	 \thanks{
	  Work supported in part by
	  ARC Discovery Project Grant
	  }\\
	  School of Mathematics\\
	  University of New South Wales\\
	  Sydney Australia \\
	  \texttt{danielc@unsw.edu.au}
	\and Colin Ingalls 
	 \thanks{
	  Work supported in part by
	  NSERC Discovery Grant}	  
       \\ 
	  Department of Mathematics and Statistics\\
	  University of New Brunswick\\
	  Fredericton, Canada\\
	  \texttt{cingalls@unb.ca}
	 }
	
%\author{Daniel Chan \\
%University of New South Wales\\
%Sydney, Australia\\
%\\
 %and \\
%Colin Ingalls
%University of New Brunswick \\
%Fredericton, Canada\\
%cingalls@unb.ca}

\date{}

\maketitle

\begin{abstract}  We discuss natural connections between three objects: 
quadratic forms with values 
in line bundles, conic bundles and quaternion orders.  We use the even Clifford
algebra \cite{Knus}, and the Brauer-Severi Variety, and other constructions
to give natural bijections between these objects under appropriate hypothesis.
We then restrict to a surface base and we express the second Chern class
of the order in terms $K^3$ and other invariants of the corresponding conic
bundle.  We find the conic bundles corresponding to minimal del Pezzo quaterion
orders and we discuss problems concerning their moduli.
\end{abstract}
\maketitle
%\tableofcontents

\begin{section}{Introduction}

In this paper, we work over a field $k$ of characteristic not equal to 2. When we speak of varieties, we mean quasi-projective varieties over the field $k$ which we assume is algebraically closed. All schemes by default will be noetherian. 

Classically, Clifford algebras over a field provide a nice construction of central simple algebras of dimension $n^2$ where $n$ is a power of two. One of our main aims is to explicitly construct quaternion orders. These are sheaves of algebras over a smooth variety say $Z$, which are locally free of rank 4 and are generically central simple over the function field $k(Z)$. A natural approach is to extend the theory of Clifford algebras to the scheme setting.  This is done
in \cite{Knus} and we apply their construction to give natural
connections between various objects. 

To motivate the scheme-theoretic generalization, recall the well known fact that terminal quaternion orders on a smooth surface $Z$ correspond to standard conic bundles on $Z$ \cite{AM}, \cite{Sarkisov}. Now a conic bundle $X$ can be written down explicitly since the relative anti-canonical embedding shows they embed in a $\P^2$-bundle, say $\P(V^*)$ for some rank 3 vector bundle $V$ on $Z$ and furthermore, $X$ is carved out by some quadratic form $Q:\Sym^2 V \rtar \L$ for some line bundle $\L$ on $Z$. It seems natural that one should be able to construct the quaternion order corresponding to $X$ using the data of this quadratic form. Now when $\L = \O_Z$ one can construct the usual Clifford algebra as the quotient of the tensor algebra $T(V)/I$ where $I$ is the ideal generated by $vw+wv - 2Q(v,w)$ for all sections $v,w \in V$. Unfortunately, this is not possible if $\L\neq \O$ but what is surprising is that the even part of the Clifford algebra $Cl_0(Q)=\O_Z \oplus \wedge^2V \otimes \L^*$ is a well-defined algebra.  We
recall this construction due to \cite{Knus} in (see \S\ref{sevencliff}). The key reason why this works is because for $v, w \in V, l \in \L^*$ we still have the perfectly legal skew commutation relation $vwl = -wvl + \langle 2Q(v,w),l\rangle$.

We are interested in the relationships between the three classes of objects below:
$$
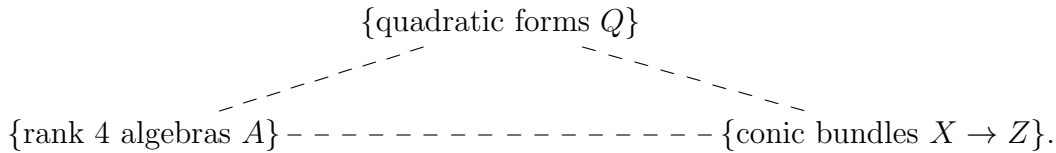
\diagram
 & \{ \text{quadratic forms}\ Q\} \dldashed \drdashed & \\  
\{ \text{rank 4 algebras}\ A \} \rdashed & & 
\{ \text{conic bundles}\  X \rtar Z\}. \ldashed
\enddiagram$$
More specifically, we consider the following questions. Which rank 4 algebras arise as even Clifford algebras? For rank 4 algebras $A$ which arise as quaternion orders, we obtain a conic bundle $BS(A)$ by taking the Brauer-Severi variety. Is this compatible with the maps $Q \mapsto Cl_0(Q)$ and $Q \mapsto X:= V(Q = 0)$ above? Given an appropriate rank 4 algebra, is there a quadratic form associated to it? Finally, given a conic bundle, can one associate a rank 4 algebra which recovers the quaternion order from its Brauer-Severi variety? After restricting the 3 classes with appropriate adjectives, we give a reasonable answer to these questions. 

We also show how natural invariants of conic bundles match those of quaternion orders.  In particular we show how the second Chern class of a quaterion order
on a surface can be expressed in terms of $-K^3$ of its Brauer-Severi
variety conic bundle and other natural invariants.

We use the correspondences above to study in particular quaternion orders which are minimal del Pezzo. These are orders which arise in the Mori program for classifying orders on surfaces \cite{CI}. They are the non-commutative analogues of del Pezzo surfaces and so deserve special attention. They have been classified using the Artin-Mumford sequence in \'etale cohomology, which can be used to show orders with prescribed ramification data exist, but give no hint as to what they look like. Now we are finally in a position to write these orders down explicitly as even Clifford algebras. Furthermore, we identify their Brauer-Severi varieties with well-known threefold conic bundles. 

We also show how natural invariants of conic bundles match those of quaternion orders.  In particular we show how the second Chern class of a quaterion order
on a surface can be expressed in terms of $-K^3$ of its Brauer-Severi
variety conic bundle and other natural invariants.

We also discuss several problems and connections between the moduli spaces of del Pezzo quaterion orders and there 
corresponding conic bundles.

The outline of the paper is as follows. 
In \S\ref{scbundle}, we review some facts about conic bundles. In particular, we recall that conic bundles on $Z$ are in bijective correspondence with orbits of ``nice'' quadratic forms under the action of $\Pic Z$. In \S\ref{sevencliff}, we recall the construction of the even Clifford algebra $Cl_0(Q)$ associated to a quadratic form $Q$ and prove some basic properties about it. 
In \S\ref{squat} we show that even Clifford algebras on rank 3 bundles have a trace function and are Cayley-Hamilton of degree two. This almost characterizes the even Clifford algebras algebraically. 
In \S\ref{squatiscliff}, we show how to recover the quadratic $Q$ from the algebraic structure of the even Clifford algebra $Cl_0(Q)$, at least under some mild additional assumptions. The key is to use the trace pairing.
In \S\ref{sbsvofcliff} we show that $Cl_0(Q)$ is the ``right'' algebra to associate to $Q$ in the sense that its Brauer-Severi variety is the conic bundle determined by $Q$. We know that forming the Brauer-Severi variety of a quaternion order is a correspondence between quaternion orders and conic bundles.  We show
an explicit inverse correspondence.  Some of the material in this paper has ``folklore status''. 

%\marginpar{Colin to rewrite rest}
In section \S\ref{sc2} we give a relation between the second Chern
class of a quaternion order and $-K^3$ of the associated conic bundle.
The rest of the paper \S\ref{sdPO} looks in depth at the case of del Pezzo and ruled quaternion orders. The del Pezzo condition depends only on ramification data and the possibilities were classified in \cite{CK,CI,AdJ}. The first task is to associate to such ramification data an appropriate quadratic form $Q$. When the centre $Z$ of the order is $\P^2$, as is the case when it is minimal del Pezzo, we may use Catanese theory \cite{Cat} with line bundle resolutions
to generate $Q$. 
%We remind the reader of the appropriate Catanese theory in the appendices. 
We compute, natural quadratic forms $Q$ associated to the ramification data of minimal del Pezzo orders and describe the corresponding Clifford algebras.
  Our theorem shows that the Brauer-Severi varieties of these Clifford algebras are just the associated conic bundles.  We identify these conic bundles with well-known descriptions of Fano three-folds described in the literature.  We
give several problems concerning moduli of these orders and their derived
categories.

\end{section}

\begin{section}{Conic Bundles}  \label{scbundle}  %\marginpar{scbundle}

In this section, we remind the reader about some basic facts concerning conic bundles and quadratic forms. 

Let $Z$ be a scheme and let $V$ be a rank $n$ vector bundle on $Z$. We will mainly be interested in the case where $Z$ is a smooth variety. 
We write $\P(V^*)$ for the scheme parametrizing rank one quotients of $V^*$.
Let $\pi: \P(V^*) \rightarrow Z$ be the projection and 
let $\O_{\P(V^*)}(H)$ be the
universal line bundle on $\P(V^*)$ associated with rank one quotients of $V^*$.
Recall that there is a canonical exact sequence
\begin{equation} \label{Pmap} 0 \rightarrow  \Omega_{\P(V^*)/Z}^1(H) \rightarrow 
\pi^* V^* \rightarrow \O_{\P(V^*)}(H) \rightarrow 0 \end{equation}
and that 
$$ \omega_{\P(V^*)/Z} \simeq \O_{\P(V^*)}(-nH) \otimes \pi^* \det V^* $$
as in \cite{H} chapter~III, Ex. 8.4.

We now introduce the notion of a quadratic form $Q$ on $V$ with values in a line bundle $\L$ on $Z$.
This is just a map $Q: \Sym^2 V \rtar \L$ so that we may view 
 $$Q \in  H^0(Z,\Sym^2 V^* \otimes \L)=
H^0(\P(V^*),\O_{\P(V^*)}(2H) \otimes \pi^*\L).$$ 
Let $X = X(Q)$ be the subscheme in $\P(V^*)$ cut out by $Q$. When $Q$ is non-zero, this is the {\it quadric bundle} associated to $Q$. When $V$ has rank 3, we shall call $X(Q)$ a {\it conic bundle}. This is more general than some definitions of conic bundles in the literature. 

The adjunction formula gives 
$$ \omega_{X/Z} \simeq \O_X((-n+2)H) \otimes \pi^*(\det V^* \otimes \L).$$

Recall there is a surjective ``symmetrizer map'' 
$$V \otimes V \rtar \Sym^2 V: v \otimes w \mapsto \frac{1}{2}( v \otimes w + w \otimes v) .$$
So sometimes, we will refer to quadratic forms $Q: V \otimes V \rtar \L$, by which we just mean one which factors through this symmetrizer map. 

Let $\widetilde{Q} : V \otimes \L^* \rightarrow V^*$ be the natural map
given by contracting with $Q$. This is the {\it symmetric matrix} associated to $Q$ which is symmetric in the sense that $\widetilde{Q}^* = \widetilde{Q} \otimes\L^*$. Note conversely that such a symmetric matrix $\widetilde{Q}$ determines a map $Q:\Sym^2 V \rtar \L$. 
When $Z$ is integral, the {\it rank} of $Q$ is just the generic rank of $\widetilde{Q}$. We say $Q$ is {\it non-degenerate} when $\widetilde{Q}$ is injective, that is, has full rank. If $Q$ is surjective, then we say that it is {\it nowhere zero}. 

There is a natural action of $\Pic Z$ on the set of quadratic forms.  Let $\M$ be 
a line bundle on $Z$. Then we obtain a new symmetric matrix 
$$ \widetilde{Q} \otimes \M: V \otimes \L^* \otimes \M \simeq 
                (V \otimes \M^*) \otimes \L^* \otimes \M^{\otimes 2}  
              \rtar (V \otimes \M^*)^* $$
which corresponds to a quadratic from on $V \otimes \M^*$ with values 
in the line bundle $\L \otimes \M^{\otimes -2}$. Note that the quadric bundle associated to $\widetilde{Q} \otimes \M$ is essentially the same as that of $\widetilde{Q}$ since they both represent essentially the same element of $H^0(Z,\Sym^2 V^* \otimes \L)$.

It turns out that if the rank $n$ of $V$ is odd then we can choose 
$\M$ to normalize $\widetilde{Q}$ as follows. Suppose for the rest of this section that $Z$ is a smooth variety. We pick a divisor $D \in \text{Div}\ Z$ with 
$$ \O(D) \simeq (\det V)^{-2} \otimes \L^n .$$
When $Q$ is non-degenerate, we can choose $D$ to be the effective 
divisor $\det \widetilde{Q} = 0$ which we note is unchanged if we alter $Q$ by a line bundle 
$\M$. 

Now $\widetilde{Q} \otimes (\det V) \otimes \L^{- \frac{n-1}{2}}$ is the map 

\diagram 
V \otimes \L^* \otimes (\det V) \otimes \L^{-\frac{n-1}{2}} \rto \dto^{\wr} 
& V^* \otimes \det V \otimes \L^{-\frac{n-1}{2}} \dto^{\wr}  \\  
(V \otimes (\det V)^* \otimes \L^{\frac{n-1}{2}}) \otimes \O(-D) \rto
& (V \otimes (\det V)^* \otimes \L^{\frac{n-1}{2}})^*    
\enddiagram

In other words, if we replace $V$ with $\widetilde{V} = V \otimes (\det V)^* \otimes \L^{\frac{n-1}{2}}$ then $\L$ gets replaced with $\O(D)$. 
%In the non-degenerate case, the normalization satisfies
%$$  (\det V)^2 \simeq \O((n-1)D)    .$$
This normalization is natural in two respects. Firstly, in the case of conic bundles, we have 
$$ \det \widetilde{V} = \det V \otimes (\det V)^{-3} \otimes \L^3 \simeq \O(D) .$$
Hence, on normalizing we may assume that $\L = \det V$ and the formula above for the relative anti-canonical bundle shows that $\omega_{X/Z}^{-1} \simeq \O_X(H)$. Secondly, the sheaf $\cok (\widetilde{Q}: V \otimes \O(-D) \rtar V^*)$ is the 2-torsion line bundle on $D$ which defines the ramification of the order. 
%(see Appendix: proposition~\ref{pcatram}). 

A conic bundle $X(Q) \rtar Z$ which is flat can be characterized intrinsically as follows. 
\begin{proposition} \label{gorsch} %\marginpar{gorsch}
Let $X$ be a Gorenstein scheme over a smooth variety $Z$ such that the fibres of $\pi:X \rtar Z$ are all (possibly degenerate) conics in $\P^2$. Then $X$ is a flat conic bundle. 
\end{proposition}
\textbf{Remark:} The converse is clear since conic bundles are hypersurfaces. In fact, one sees easily that flat conic bundles are precisely those of the form $X(Q)$ where $Q$ is nowhere zero. 

\begin{pf}
Note  $\pi$ is flat since $Z$ is smooth, $X$ is Gorenstein and the fibres of $\pi$ are all 1-dimensional. Also, the relative anti-canonical bundle $\omega_{X/Z} := \omega_X \otimes \pi^*\omega_Z^{-1}$ is flat over $Z$. Grauert's theorem and the condition on the fibres now ensure $V^*:=\pi_*\omega_{X/Z}^{-1}$ is a vector bundle of rank 3 and we have a relative anti-canonical embedding $X \hookrightarrow \P_Z(V^*)$. Computing fibre-wise, we see that the corresponding line bundle $\O_{\P(V^*)}(X) \simeq \O_{\P(V^*)}(2H) \otimes \pi^* \L$ for some $\L \in \Pic Z$. Now $X$ is given by a section of this bundle so determines up to scalar a quadratic form $Q= Q(X) \in \Hom_Z(\Sym^2 V, \L)$. 
\end{pf}

The argument above shows that for flat conic bundles we have $\pi_* \omega_{X/Z}^{-1}$ is a rank three vector bundle. This is true in general. 
\begin{lemma}  \label{lrantican}  %\marginpar{lrantican}
Let $Q$ be a quadratic form on a rank three vector bundle $V$ with associated conic bundle $X$. Then 
$$\pi_* \omega_{X/Z}^{-1} \simeq V^* \otimes \det V \otimes \L^*  .$$
In particular if $Q$ is normalized, then $\pi_* \omega_{X/Z}^{-1} \simeq V^*$.
\end{lemma}
\begin{pf}
Consider the exact sequence of sheaves on $\P(V^*)$
$$ 0 \rtar \O_{\P(V^*)/Z}(-H) \otimes \pi^* \L^* \rtar \O_{\P(V^*)/Z}(H) \rtar \O_{X/Z}(H) \rtar  0.$$
For $i=0,1$ we have 
$$ R^i\pi_* (\O_{\P(V^*)/Z}(-H) \otimes \pi^* \L^*) = R^i\pi_* \O_{\P(V^*)/Z}(-H) \otimes \L^* =0 $$ 
so the long exact sequence in cohomology gives $\pi_* \O_{X/Z}(H) = \pi_* \O_{\P(V^*)/Z}(H) = V^*$. 
The adjunction formula above gives for $X=X(Q)$ 
$$\pi_* \omega_{X/Z}^{-1} \simeq \pi_* (\O_{X/Z}(H) \otimes \pi^* \det V \otimes \pi^* \L^*) \simeq 
V^* \otimes \det V \otimes \L^*  .$$
\end{pf}

%\vspace{2mm}

%We get a natural map \marginpar{Colin to complete}
%$$\O_{\P(V^*)}(-H) \otimes \pi^*\L^* \rightarrow 
%\pi^*(V \otimes \L^*)$$ 
%by dualizing \ref{Pmap} and twisting by $\pi^*L$.
%We compose this map with 
%$$\pi^*\widetilde{Q} : \pi^*(V \otimes \L^*) \rightarrow \pi^*V^*$$
%and call the cokernel of the composition $I$.
% 
%We obtain a commutative diagram
%$$ 
%\begin{CD}
%\O_{\P(V^*)}(-H) \otimes \pi^*\L^* &  @= &   \O_{\P(V^*)}(-H) \otimes \pi^*\L^* \\
%@VVV &   &                @VVV \\
%\pi^*(V \otimes \L^*) & @>>> & \pi^*V^* & @>>> & \pi^*\cok \widetilde{Q}  \\
%@VVV &   & @VVV & & @| \\
%T_{\P(V^*)/Z}(-H) \otimes \pi^*L^* &  @>>>  & I & @>>> & \pi^* \cok \widetilde{Q}
%\end{CD}
%$$
%where all the above sequences are short exact with the zeroes ommitted.

%Restricting
%the above sequence to $X$ gives the diagram
%$$ 
%\begin{CD}
%\O_X(-H) \otimes \pi^*\L^* &  @= &   \O_X(-H) \otimes \pi^*\L^*) \\
%@VVV &   &                @VVV \\
%\Omega_{\P(V^*)/Z}(H)|_X & @>>> & \pi^*V^*|_X & @>>> & \O_X(H) \\
%@VVV &   & @VVV & & @| \\
%\Omega_{X/Z}(H) &  @>>>  & I|_X & @>>> & \O_X(H) 
%\end{CD}
%.$$

%There are several natural bundles on conic bundles
%that we will need.  There is the 

\end{section}

\begin{section}{Even Clifford Algebras}  \label{sevencliff} %\marginpar{sevencliff}
We now recall the construction of the even Clifford algebra of a quadratic
form with values in a line bundle due to Bichsel and Knus \cite{Knus}.  Let $Z$ be a scheme and $Q:\Sym^2 V \rtar \L$ be a quadratic form on a rank $n$ vector bundle $V$ with values in the line bundle $\L \in \Pic Z$. 
When $\L = \O$, there is the well-known construction of the Clifford algebra, which is a sheaf of $\Z/2\Z$-graded $\O_Z$-algebras of rank $2^n$.  A version of the even part of the Clifford algebra can be defined.  There is also
a version of the odd part of the Clifford algebra that is a module over
the even part, but the even
and odd parts do not form an algebra.  To construct the even part we 
proceed as follows:b

We first consider two $\Z$-graded $\O_Z$-algebras: the tensor algebra $T(V) = \oplus T(V)_i$ and $\oplus_{j \in \Z} \L^j$. Tensoring these two algebras together gives a bigraded algebra 
$$  T(V,\L) := T(V) \otimes_Z (\oplus_{j \in \Z} \L^j)  .$$
Now $\Sym^2 V \subset T(V)_2 = V \otimes_Z V$ so we may consider $Q$ as a relation in $T(V,\L)$ and define the {\it total Clifford algebra} $Cl_{\bullet}(Q)$ to be the quotient of $T(V,\L)$ with defining relation $Q$. More precisely, let $I\triangleleft T(V,\L)$ be the two-sided ideal generated by sections of the form 
$t - Q(t)$ for all $t \in \Sym^2 V$. Then 
$$ Cl_{\bullet}(Q):= T(V,\L)/ I   .$$
If $\widetilde{Q}$ is the symmetric matrix associated to $Q$ then we also write $Cl_{\bullet}(\widetilde{Q})$ for $Cl_{\bullet}(Q)$. 

Of course, $Cl_{\bullet}(Q)$ is no longer bigraded. However, if we give $V$ degree 1 and $\L$ degree 2, then the relation is homogeneous of degree 2 so $Cl_{\bullet}(Q)$ is $\Z$-graded. The degree zero part $Cl_0(Q)$ is called the {\it even Clifford algebra} since, when $\L = \O$, it is the even part of the usual Clifford algebra. 

Recall from section~\ref{scbundle} that $\Pic Z$ acts on $Q$. Though altering $Q$ by a line bundle $\mathcal{M} \in \Pic Z$ affects $Cl_{\bullet}(Q)$, it does not affect $Cl_0(Q)$. 

We need a result concerning the classical Clifford algebra of a quadratic form $Q: V \otimes V \rtar \O_Z$ which defined by $Cl(Q) = T(V)/I$ where $I$ is the ideal generated by sections $t- Q(t)$ for $t \in \Sym^2 V$.

\begin{proposition}  \label{pClstar}  %\marginpar{pClstar}
Let $V = \O_Z^2$ and $Q:V \otimes V \rtar \O_Z$ be a quadratic form. Then $Cl(Q)^* \simeq Cl(Q)$ as left and right $Cl(Q)$-modules.
\end{proposition}
\begin{pf}
Write $V = \O_Z x \oplus \O_Z y$ and note that 
$$Cl(Q) = \O_Z \oplus \O_Z x \oplus \O_Z y \oplus \O_Z xy.$$
Let $\xi: Cl(Q) \rtar \O_Z$ be projection onto $\O_Z xy$. Then one computes 
readily that the left and right modules generated by $\xi$ are both the whole of $Cl(Q)^*$. 
\end{pf}

\begin{proposition}  \label{peveniscliff}  %\marginpar{peveniscliff}
Let $Z$ be the spectrum of a local ring with closed point $p$ and $Q:V\otimes V \rightarrow \O_Z$
be a quadratic form on a rank $n$ vector bundle $V$. Suppose the induced quadratic form $Q\otimes_Z k(p): V_p \otimes_k V_p \rtar k(p)$ is non-zero where $V_p = V \otimes_Z k(p)$.  Then there is a rank $n-1$ sub-bundle $V'<V$ and a quadratic form $Q':V' \otimes V' \rtar \O_Z$ of rank $\rank Q - 1$ such that $\Cl_0(Q) \simeq Cl(Q')$. 
\end{proposition}
\begin{pf}
Since $Q_p$ does not vanish, we can find a section $u \in V$ such that $Q(u,u)\in \O_Z^*$.  The map $Q(-,u)$ from $V \rightarrow \O$ is surjective since $Q(u,u)$ is a unit.  So the kernel $V'$ is locally
free of rank $n-1$. Now $Cl_0(Q)$ is generated by $uV'$ and, identifying $uV'$ with $V'$ in the natural way, we see that $Cl_0(Q) \simeq Cl(Q')$ where $Q' = -Q(u,u) Q|_{V'}$. 
\end{pf}

The Azumaya locus is given by the following proposition.

\begin{proposition}  \label{pazulocus}   %\marginpar{pazulocus}
Suppose $Q: V \otimes \L^* \rtar V^*$ is a symmetric matrix and $V$ is odd dimensional. 
The Azumaya locus of the even Clifford algebra $Cl_0(Q)$ is the non-degeneracy locus of $Q$, that is, the open set $\det Q \neq 0$.
\end{proposition}
\begin{pf}
The Azumaya locus is where the fibres are central simple algebras. Hence we may restrict to a point so that $\O_Z$ is a field. It is known classically that the even Clifford algebra in this case is central simple if and only if $Q$ is non-degenerate.
\end{pf}

To study the even and the total Clifford algebra, we notice that $Cl_{\bullet}(Q)$ also has an ascending filtration where the $i$-th filtered piece is
$$ F^i Cl_{\bullet}(Q) = \mbox{im}: (\oplus_{l \leq i} T(V)_l) \otimes (\oplus_{j \in \Z} \L^j) \rtar Cl_{\bullet}(Q)  .$$
The associated graded algebra is then easily seen to be
$$ gr Cl_{\bullet}(Q) = \wedge^{\bullet} V \otimes_Z (\oplus_{j \in \Z} \L^j)  .$$
This shows in particular that $Cl_0(Q)$ is locally free of rank $2^{n-1}$. 

Recall that the wedge product induces a perfect pairing $\wedge^r V \otimes \wedge^{n-r} \rtar \det V$ so if $n$ is odd, there is a duality between $\wedge^{\text{even}} V$ and $\wedge^{\text{odd}} V$. We will see shortly that the same is true for Clifford algebras. 

For the rest of this section, we assume that $n=3$ which corresponds to conic bundles and algebras of rank 4. In this case, the filtration gives an exact sequence of sheaves
$$ 0 \rtar  \O_Z \rtar  Cl_0(Q) \rtar \wedge^2 V \otimes \L^{-1} \rtar 0 .$$
If $Q$ is non-degenerate, then generically $Cl_0(Q)$ is central simple so the even Clifford algebra is an order. It follows that the reduced trace splits the above sequence and, writing $sCl_0(Q)$ for the traceless part of $Cl_0(Q)$, we have 
$$ Cl_0(Q) = \O_Z \oplus s Cl_0(Q) , \hspace{1cm} sCl_0(Q) \simeq \wedge^2 V \otimes \L^{-1} .$$
When $Q$ is normalized so that $\L = \det V$, we further have $sCl_0(Q) \simeq V^*$. 

The next result shows how to recover the total Clifford algebra from the even part. 
\begin{proposition}  \label{ptotalcliff}  %\marginpar{ptotalcliff}
Consider a total Clifford algebra $Cl_{\bullet}(Q)$ where $Q$ is a normalized quadratic form on a rank 3 vector bundle and let $A=Cl_0(Q)$ be the even Clifford algebra. The graded decomposition of $Cl_{\bullet}(Q)$ can be rewritten as
$$  Cl_{\bullet}(Q) = (\bigoplus_{j \in \Z} A \otimes_Z \L^j) \oplus (\bigoplus_{j \in \Z} A^* \otimes_Z \L^j)$$
where $A^*$ sits in degree 1. The decomposition is as of $A$-modules. Moreover, in this description $(A/\O_Z)^* < A^*$ corresponds to $V < Cl_1(Q)$. 
\end{proposition}
\begin{pf}
The filtration on the third graded component of $Cl_{\bullet}(Q)$ gives the exact sequence 
$$ 0 \rtar V \otimes \L \rtar Cl_3(Q) \rtar \wedge^3 V \rtar 0  .$$
Since $Q$ is normalized, we may identify $\L$ with $\wedge^3 V$, to see that multiplication in the total Clifford algebra gives a pairing 
$$ Cl_1(Q) \otimes_A Cl_2(Q) \rtar Cl_3(Q) \rtar \L .$$
It is a perfect pairing since it is compatible with the perfect pairing on $gr Cl_{\bullet}(Q)$. 
Tensoring with $\L^{-1}$ shows that $Cl_1(Q) = Cl_0(Q)^*$. It also shows that $(A/\O_Z)^* < A^*$ corresponds to $V$. 
\end{pf}

The following result will be useful in the next section.

\begin{lemma}  \label{lcl1}  %\marginpar{lcl1}
Let $V$ be a rank three vector bundle and $Q: V \otimes V \rtar \L$ a quadratic form with values in a line bundle $\L$. Suppose that $Q$ is non-degenerate and nowhere zero. Writing $A$ for the even Clifford algebra, we have an $A$-bimodule isomorphism $A^* \otimes_A A^* \simeq A \otimes \L$. Furthermore, the isomorphism maps $\Sym^2 (A/\O_Z)^*$ onto $\O_Z \otimes \L$.
\end{lemma}
\begin{pf}
Consider the bimodule morphism given by multiplication in the total Clifford algebra 
$$  \mu: A^* \otimes_A A^* = Cl_1(Q) \otimes_A Cl_1(Q) \rtar Cl_2(Q) = A \otimes \L .$$
Note that $Cl_2(Q) = V^2 + \L$ where $V^2$ denotes sums of products of elements in $V$. Since $Q$ is nowhere zero, we in fact have $V^2 \supset \L$ so $\mu$ is clearly surjective. Note that locally, $A$ is Clifford by proposition~\ref{peveniscliff} and the assumption that $Q$ is nowhere zero. So locally on $Z$, proposition~\ref{pClstar} shows that $A^* \simeq A$ as a left and right $A$-module. Hence $A^* \otimes_A A^*$ is locally isomorphic to $A$ and $\mu$ induces the desired isomorphism. 
\end{pf}

\textbf{Remark:} If $A$ is a maximal order in a 4-dimensional central simple algebra then ramification theory says that $A^* \otimes_A A^* \simeq A(D)$ where $D$ is the discriminant. Presumably the above map must surely coincide with this one.

\end{section}

\begin{section}{Quaternion Algebras}  \label{squat}  %\marginpar{squat}

Not all locally free algebras of rank four occur as even Clifford algebras. We give a partial intrinsic characterization of these algebras. 
Let $Z$ be a scheme.  Usually, $Z$ will be integral and $A$ will be an $\O_Z$-algebra that is locally free of rank four. 

\begin{definition}  \label{dquat}   %\marginpar{dquat}
We say that an $\O_Z$-algebra $A$ is {\it quaternion} if it is locally free of rank four and there is a linear trace function $\mbox{tr}: A \rtar \O_Z$ such that 
\begin{enumerate}
\item $\frac{1}{2}\mbox{tr}$ splits the natural inclusion $\O_Z \rtar A$. 
\item any section $a \in A$ satisfies a quadratic relation of the form $a^2 - tr(a) a + g=0$ where $g \in \O_Z$.
\end{enumerate}
\end{definition}
\begin{definition}
Suppose that $Z$ is a normal integral scheme. A {\it quaternion order} is an $\O_Z$-algebra $A$ that is locally free of rank four and such that $k(Z) \otimes_Z A$ is a central simple $k(Z)$-algebra. The definition is justified by the fact that the reduced trace function satisfies conditions 1) and 2) so $A$ is quaternion. 
\end{definition}
\textbf{Remark:} The conditions 1) and 2) above define the trace uniquely  when $Z$ is integral (see for example the proof of the proposition below). Furthermore, in 2) we have $g = \frac{1}{2}( (\text{tr}\ a)^2 - \text{tr}\ a^2)$.

Here are some basic facts. For the definition of Cayley-Hamilton algebras, see [leB, \S1.6]. 
\begin{proposition}  \label{paltquat}  %\marginpar{paltquat}
Let $A$ be an $\O_Z$-algebra that is locally free of rank four. Then it is quaternion if and only if we can split the sheaf $A = \O_Z \oplus sA$ such that for every $x \in sA$ we have $x^2 \in \O_Z$. If $A$ is a quaternion algebra on an integral scheme $Z$, then the trace pairing $A \times A \rtar \O_Z: (x,y) \mapsto \mbox{tr}(xy)$ is symmetric so $A$ is Cayley-Hamilton of degree two. 
\end{proposition}
\begin{pf}
Assume $A$ has a splitting $A = \O_Z \oplus sA$ as above. The splitting defines a trace map which satisfies the conditions of definition~\ref{dquat} since we are assuming $x^2 \in \O_Z$ for every $x \in sA$. To prove symmetry of the trace pairing it suffices, since $\O_Z$ is central, to show $\mbox{tr}(xy) = \mbox{tr}(yx)$ for all linearly independent $x,y \in sA$. Let $t = xy+yx$ which lies in $\O_Z$ since 
$$ t = (x+y)^2 - x^2 - y^2 \in \O_Z .$$
Consider first the case where $xy \notin \O_Z$ so its trace is determined by condition~2) of definition~\ref{dquat}. We have $(xy)^2 = -x^2y^2 + txy$. Since $x^2,y^2 \in \O_Z$ we see $\mbox{tr}(xy) = t$ and a symmetric computation shows $\mbox{tr}(yx) = t$. Suppose on the other hand that $xy \in \O_Z$ so also $yx \in \O_Z$. Then $(x^2)y + x(yx) = xt$ so linear independence of $x,y$ force $x^2 = 0,yx = t$. A symmetric computation shows $xy = t$ so in fact we will have $xy=yx=0$. 
\end{pf}

The next result immediately suggests a strong relationship between quaternion algebras and Clifford algebras. 
\begin{proposition}  \label{ploccliff}  %\marginpar{ploccliff}
Let $Z$ be the spectrum of a local ring with closed point $p$, and  $A$ be a quaternion $\O_Z$-algebra. Suppose that $A\otimes_Z k(p)$ is generated as a $k(p)$-algebra by two elements $x,y\in A$. Then $A= Cl(Q)$ for some quadratic form on $V = \O_Z (x - \frac{1}{2}tr x) \oplus \O_Z (y - \frac{1}{2}tr y)$. 
\end{proposition}
\begin{pf}
Replacing $x,y$ with $x - \frac{1}{2}tr x,y - \frac{1}{2}tr y$, we may suppose that $x,y \in sA$. Hence $a:=x^2,b:=y^2,c:=\frac{1}{2}(xy+yx) \in \O_Z$. Now $x,y$ generate $A \otimes_Z k(p)$ so $1,x,y,xy$ must form a $k(p)$-basis. They are thus also an $\O_Z$-basis for $A$. If we set $V = \O_Z x \oplus \O_Z y$, then $A = Cl(Q)$ where 
$$ Q = \begin{pmatrix} a & c \\ c & b \end{pmatrix}.$$
\end{pf}

Restricting a quaternion algebra to a closed subscheme $Y \subset Z$ gives a quaternion algebra on $Y$. The closed fibres of a quaternion algebra are thus prescribed by the following result.

\begin{theorem}  \label{tclassquat}  %\marginpar{tclassquat}
Let $A$ be a quaternion $k$-algebra. Then $A$ is isomorphic to one of the following algebras:
\begin{enumerate}
\item a central simple Clifford algebra $k\langle x,y\rangle/(x^2-a,y^2-b,xy+yx)$ for some $a,b \in k^*$.
\item a Clifford algebra of form $k\langle x,y\rangle/(x^2-a,y^2,xy+yx)$ for some $a \in k^*$.
\item the Clifford algebra $k\langle x,y\rangle/(x^2,y^2,xy+yx)$.
\item the commutative algebra $k[x,y,z]/(x,y,z)^2$.
\item the quiver algebra of the Kronecker quiver with two arrows: 
%$k\langle e,x,y\rangle/(e^2-e,ex-x,ey-y,ye,xe,y^2,x^2,xy,yx)$.
$$\begin{pmatrix}
k & V \\ 
0 & k
\end{pmatrix}$$
where $V$ is a two dimensional vector space over $k$. 
\end{enumerate}
The only Clifford algebras amongst these are 1),2) and 3). The algebras 1),2),3),4) are the even Clifford algebras $Cl_0(Q)$ where $Q$ has rank 3,2,1,0 respectively. 
\end{theorem}
\begin{pf} We will use proposition~\ref{ploccliff} repeatedly without further comment. Consider first the trace pairing $P$ on the traceless part $sA$ of $A$. If $P$ has rank $\geq 2$, then we can find an orthogonal subset $\{x,y\} \subset sA$ with $a:=x^2,b:=y^2 \in k^*$. Now $xy \in sA$ so $xy + yx \in \O_Z$ must be zero. Thus $(xy)^2 \neq 0$ and $xy$ is orthogonal to $kx + ky$. It follows that $\{1,x,y,xy\}$ is a basis for $A$ from which we see it is 2-generated and so must be the central simple Clifford algebra 1). If $\rank P = 1$ then there are several cases. We pick as before $x \in sA$ with $x^2=a \in k^*$ and let $V < sA$ be the orthogonal complement of $x$. If $a$ is not a square in $k$, then $k[x]$ is a field extension of $k$ and $A$ is generated by $x$ and any non-zero $y \in V$ so we are in case 2). Suppose now that $a \in k^{*2}$ so on scaling $x$ we may suppose that $x^2=1$. Note then that $k[x] \simeq k(\Z/2\Z) \simeq k \times k$ and that $V$ is a $k[x]$-bimodule. If there is an isomorphism of left modules $V \simeq k[x]$, then we may pick $y \in V$ to be a non-zero non-eigenvector with respect to $x$. It follows that $\{1,x,y,xy\}$ is a basis for $A$ so we are in case 2). In the other case, $V$ is a single eigenspace for $x$, so using the idempotents in $k[x]$ we obtain a Peirce decomposition for $A$ which gives the algebra in 5). Finally, if $\rank P = 0$ then for any basis $\{x,y,z\} \subset sA$ we have $0=x^2=y^2=z^2 = xy+yx=yz+zy=xz+zx$. If all the products $xy,zx,yz$ are zero, then we are in case 4) so suppose, without loss of generality that $xy\neq 0$. To show that case 3) holds, it suffices to prove that $xy$ is linearly independent from $x,y$. But if $xy = ax + by$ then the equalities $0=x^2y=xy^2$ force $a=b=0$ so $A$ is indeed the Clifford algebra of 3).

One can calculate the even Clifford algebra $Cl_0(Q)$ when $Q$ has rank 3,2,1,0 to obtain the algebras in 1),2),3),4) respectively. Hence 5) is not an even Clifford algebra. 
\end{pf}

\textbf{Remark:} a) The algebra in 5) above does not occur in the case of orders for the following reason. Suppose that $Z$ is the spectrum of a complete local ring with closed point $p$ and residue field $k(p)= k$. If $A\otimes_Z k(p)$ is the quaternion algebra in 5), then it has a non-trivial idempotent $e$ such that $e(A \otimes_Z k(p))(1-e) = 0$. Now idempotents may be lifted to $Z$ to give a Peirce decomposition
$$A \simeq 
\begin{pmatrix}
\O_Z & V \\
0 & \O_Z
\end{pmatrix}$$
where $V$ is a vector bundle of rank two on $Z$. In other words, one cannot deform the quiver algebra of the Kronecker quiver into the matrix algebra.

b) Any Clifford algebra formed from a symmetric matrix $\widetilde{Q}: V \rtar V^*$ is isomorphic to the even Clifford algebra formed from the symmetric matrix $\widetilde{Q} \oplus 1: V \oplus \O \rtar V^* \oplus \O$. The theorem above shows that, even locally on $Z$, an even Clifford algebra may not be a Clifford algebra. 

%As a result of the proposition, any quaternion algebra $A$ which is generated by two elements at every closed point will be said to be {\it locally Clifford}. 

\begin{definition}
We say that an algebra $A$ is {\it locally quaternion (respectively, Clifford, or even Clifford)} if the localization of $A$ at any point is quaternion (respectively, Clifford, or even Clifford). 
\end{definition}

\begin{proposition}  \label{plocalglobal}  %\marginpar{plocalglobal}
Suppose $Z$ is an integral scheme. 
Any locally quaternion algebra is quaternion. Any locally even Clifford algebra of rank four is quaternion. 
\end{proposition}
\begin{pf}
Let $A$ be a quaternion algebra and $A = \O_Z \oplus sA$ be the splitting induced by the trace. Then 
$$ \{ a \in A| a^2 \in \O_Z \} = \O_Z \cup sA  $$
so the subbundle $sA$ is uniquely defined. In particular, any locally quaternion algebra is quaternion. 

Consider the even Clifford algebra $Cl_0(Q)$. To show it is quaternion, it suffices to work locally so we may assume $Q:V \otimes V \rtar \O_Z$ is given by a matrix $(q_{ij})$ with respect to a basis $\{x_1,x_2,x_3\}$ for $V$. Recall this means $Cl_0(Q)$ has defining relations 
$$ \frac{1}{2}(x_ix_j + x_jx_i) = q_{ij}, \hspace{1cm} \mbox{for all}\ i,j. $$
Then 
$$sA := \O_Z(x_1x_2 - q_{12}) \oplus \O_Z (x_2x_3 - q_{23}) \oplus \O_Z (x_3x_1 - q_{31})$$ 
is a complement to $\O_Z$ with which we can apply proposition~\ref{paltquat} to show that $A$ is quaternion. 
\end{pf}

We can finally give an intrinsic characterization of the even Clifford algebras of nowhere zero quadratic forms.

\begin{theorem}  \label{tiscliff}   %\marginpar{tiscliff}
Let $A$ be an $\O_Z$-algebra where $Z$ is an integral scheme. The following are equivalent.
\begin{enumerate}
\item $A$ is a quaternion algebra such that for every $p \in Z$ closed, the algebra $A \otimes_Z k(p)$ is generated by two elements. 
\item $A$ is a locally Clifford algebra of rank 4.
\item $A \simeq Cl_0(Q)$ for some nowhere zero quadratic form $Q: V \otimes V \rtar \L$ on a rank three vector bundle $V$ with values in the line bundle $\L$. 
\end{enumerate}
\end{theorem}
\begin{pf}
We assume first that 3) holds and prove 1). Proposition~\ref{peveniscliff} and the fact that $Q$ is nowhere zero implies that $Cl_0(Q)$ is locally Clifford and hence locally generated by two elements. It is quaternion by proposition~\ref{plocalglobal}. The implication 1) $\implies$ 2) is proposition~\ref{ploccliff}. 

Finally, we assume 2) and prove 3). We know that locally, $A$ is an even Clifford algebra and must show this holds globally. We construct the total Clifford algebra first. The Clifford algebra will be built from the rank three vector bundle $V = (A/\O_Z)^*< A^*$. Finding the line bundle $\L$ is more subtle. Locally on $Z$, we know by lemma~\ref{lcl1} that there is an $A$-bimodule isomorphism $A^* \otimes_A A^* \simeq A$ which maps $\Sym^2 V$ onto $\O_Z$. Now local computations show $Z(A) = \O_Z$ so bimodule isomorphisms  $A \rtar A$ are given by multiplication by  sections of $\O_Z^*$. These give transition functions which define a line bundle $\L \in \Pic Z$ allowing us to glue these isomorphisms to a global isomorphism $A^* \otimes_A A^*  \simeq A\otimes \L$. This in turn gives the ``total Clifford algebra''
$$  Cl_{\bullet}(Q) = (\bigoplus_{j \in \Z} A \otimes_Z \L^j) \oplus (\bigoplus_{j \in \Z} A^* \otimes_Z \L^j).$$
Looking locally at the isomorphism $A^* \otimes_A A^* \simeq A \otimes \L$, we see that the multiplication maps $\Sym^2 V$ into $\L$. Consequently, multiplication gives a global map $Q:\Sym^2 V \rtar \L$. It is now clear that $A \simeq Cl_0(Q)$.
\end{pf}
The theorem shows that, under the nowhere zero assumption, the property of being an even Clifford algebra is purely local on the base $Z$. This mimics the case of conic bundles.

\textbf{Remark:} If an algebra is locally Clifford in codimension one, then often one can apply the previous result on the Clifford locus and extend across codimension two points. For example, let $A$ be a maximal order of rank 4 on a smooth surface. Then $A$ is reflexive hence locally free by Auslander-Buchsbaum. Purity of the branch locus and \'etale local descriptions of $A$ at codimension one points show that $A$ is locally Clifford outside some closed subset $Y \subset Z$ of codimension at least two. Hence $A$ is isomorphic to $Cl_0(Q)$ on $Z-Y$ for some quadratic form $Q:V \otimes V \rtar \L$. Now $V,\L, Q$ all extend uniquely to $Z$ and, since $A$ is reflexive, $A$ is determined completely by its structure on $Z-Y$. Hence $A$ is globally even Clifford.

\end{section}

\begin{section}{Trace pairing of even Clifford algebras} \label{squatiscliff} %\marginpar{squatiscliff}

In the last section we saw that quaternion algebras that are locally 2-generated are even Clifford algebras. However the method of proof did not explicitly find the quadratic form to construct the Clifford algebra. In this section, following a suggestion of Johan de Jong, we find the quadratic form in terms of the trace pairing. Throughout let $Z$ be an integral scheme.

We start by computing the trace pairing of an even Clifford algebra. Consider a normalized symmetric matrix $Q:V \otimes (\det V)^* \rtar V^*$ on a rank 3 vector bundle $V$, so that $sCl_0(Q) \simeq \wedge ^2 V \otimes (\det V)^* \simeq V^*$. First note that 
$$\wedge^2 Q^*: \wedge^2 V \rtar \wedge^2 V^* \otimes (\det V)^{\otimes 2} $$
which via the cross product isomorphism induces a map 
$$ \wedge^2 Q^*: V^* \otimes \det V \rtar V \otimes \det V .$$
We define the adjoint of $Q$ to be 
$$\Adj Q:=\wedge^2 Q^* \otimes \det V^*: V^* \rtar V .$$

\begin{proposition}
The trace pairing on $sCl_0(Q)$ is given by $-\Adj Q$. 
\end{proposition}
\begin{pf}
It suffices to assume that $Z= \Spec R$ is affine and that $V = Rx \oplus Ry \oplus Rz$. We let $\{x^*,y^*,z^*\}$ be the dual basis for $V^*$ and let 
$\d:=x^* \wedge y^* \wedge z^* \in \det V^*$. We represent $Q:V \otimes (\det V)^* \rtar V^*$ by the matrix 
$$\begin{pmatrix}
 a & d & e \\
 d & b & f \\
 e & f & c
\end{pmatrix}$$
with respect to the bases $\{x\otimes\d,y\otimes\d,z\otimes \d\}$ and $\{x^*,,y^*,z^*\}$. The cross product isomorphism $\wedge^2 V \otimes \det V^* \rtar V^*$ identifies 
$$ x \wedge y \otimes \d \leftrightarrow z^*,
   z \wedge x \otimes \d \leftrightarrow y^*,
   y \wedge z \otimes \d \leftrightarrow x^* .$$
To identify these with $sCl_0$ we need to be a little careful with our notation. We view $Cl_0$ as a subquotient of the tensor algebra $R\langle x,y,z,\d \rangle$ and so write elements of $Cl_0$ as polynomials in the noncommuting variables $x,y,z,\d$. Then if $\bar{z}$ is the element of $sCl_0$ corresponding to $z^*$ we have 
$$ \bar{z} = xy\d - \frac{1}{2} \tr xy\d  .$$
We compute the trace term as follows
$$\begin{aligned}
xy\d xy\d & =  -x^2y^2\d^2 + 2Q(x,y)xy\d \\
          & =  -Q(x,x)Q(y,y)   + 2Q(x,y)xy\d \\
          & =  -ab   + 2d xy\d 
\end{aligned}$$
This gives $\tr xy\d = 2d$ and so the elements of $sCl_0$ corresponding to $z^*,y^*,x^*$ are 
$$ \bar{z} = xy\d - d, 
   \bar{y} = zx\d - e, 
   \bar{x} = yz\d - f . $$
Completing the square gives 
$$\bar{z}\bar{z} = -ab + d^2 = 
-\langle z^*,\Adj z^*\rangle \ \mbox{so}\ 
\tr \bar{z}\bar{z} = -\langle z^*,\Adj z^*\rangle. $$
Similarly, 
$$\begin{aligned}
\bar{y}\bar{z} &= (zx\d - e)(xy\d - d)  \\
               &= a zy\d -exy \d 
                 -dzx\d +ed
\end{aligned}$$
so 
$$ \tr \bar{y}\bar{z} = af - ed = 
                            -\langle y^*,\Adj z^*\rangle .$$
A similar calculation for the other pairings between basis elements gives the result.
\end{pf}

The next proposition ``solves'' for the quadratic form in terms of the trace pairing.
\begin{proposition} \label{pfindQ} 
Suppose that $\Pic Z$ has no 2-torsion. Let $A$ be a quaternion order with trace pairing $P:sA \rtar (sA)^*$. Then there is an even Clifford algebra $Cl_0(Q)$ with the same trace pairing as $A$. If $Q$ is normalized, then it is given by the formula $Q = \sqrt{-1} \Adj P (\det P)^{-\frac{1}{2}}$ which is to be interpreted as in the proof below. 
\end{proposition}
\begin{pf}
Let $D$ be the discriminant divisor of $A$.
Let $V = (sA)^*$ and $P:V^* \rtar V$ denote the trace pairing. The classical ramification theory of orders shows that $\det P:\det V^* \rtar \det V$ induces a map $\det P\otimes \det V:\O \rtar (\det V)^{\otimes 2} \simeq \O(2D)$. Hence by our assumption on $\Pic Z$, there is, up to $\pm 1$, a unique isomorphism $\det V \simeq \O(D)$ compatible with $\det P$. We write $(\det P)^{\frac{1}{2}}:\O \rtar \O(D) \simeq \det V$.

Generically, we can solve $P=-\Adj Q$ for $Q$  as follows. Note first that $\det P = - (\det Q)^2$ so 
$$ \det Q = (-\det P)^{\frac{1}{2}}.$$
Cramer's rule gives 
$$P=-Q^{-1} \det Q .$$
Hence using Cramer's rule again gives 
$$Q = -P^{-1} \det Q = \sqrt{-1} \Adj P (\det P)^{-\frac{1}{2}}.$$
The solution is unique up to $\pm 1$. 

We need to show that this solution is well-defined globally. Consider the following diagram
$$\diagram
V \otimes \det V^* \rto^{\Adj P} \drto^{Q} & V^* \otimes \det V  \\
  & V^* \uto_{V^* \otimes (\det P)^{1/2}} 
\enddiagram$$
 It suffices to show that $\Adj P$ lies in the image of $V \otimes (\det P)^{1/2}$ for then we can define $Q$ to make the above diagram commute. It suffices to check this locally in codimension one. Now $(\det P)^{\frac{1}{2}}$ is an isomorphism away from $D$ so let $R$ be a discrete valuation ring corresponding to the local ring at a prime divisor of $D$ and let $\m$ be its maximal ideal. Using Gram-Schmidt over a discrete valuation ring, we may diagonalize $P$ so that it takes the form 
$$P = \begin{pmatrix}
 a & 0 & 0 \\
 0 & b & 0 \\
 0 & 0 & c
\end{pmatrix}.$$
Now the ideal of $R$ generated by $\det P$ is $\m^2$ so, without loss of generality, we may assume that $a,b \in \m, c \notin \m$. Then 
$$\Adj P = \begin{pmatrix}
 bc & 0 & 0 \\
 0 & ac & 0 \\
 0 & 0 & ab
\end{pmatrix} \in \m R^{3 \times 3}.$$
This completes the proof of the proposition.
\end{pf}
The next lemma shows that the trace pairing gives most of the multiplicative structure of a quaternion order. 

\begin{lemma}
Let $A,A'$ be quaternion orders with the same underlying sheaf, say $\O \oplus V^*$. Assume that the trace pairing on $A,A'$ are the same and that $V^*$ is the trace zero part of $A,A'$. Then the identity morphism $A \rtar A'$ of sheaves is either an isomorphism of orders or an anti-isomorphism.
\end{lemma}
\begin{pf}
Let $K=k(Z)$. It suffices to show that generically, the identity map $\id: A \otimes K \rtar A' \otimes K$ is either an algebra isomorphism or anti-isomorphism. Pick a $K$-basis $\{x,y,z\}$ for $V^*\otimes K$ which is orthogonal with respect to the trace pairing. Note that in either $A$ or $A'$ we have 
$$  xy + yx = (x+y)^2 - x^2 - y^2  \in K$$
since $x+y,x,y$ all have trace zero. We see consequently that 
$$ xy+yx = \tr xy = 0$$
by hypothesis so $x,y,z$ all skew-commute in both $A$ and $A'$.

We first show that the multiplication on $A \otimes K$ is completely determined by the trace pairing and the scalar $\tr (xyz) \in K$. We compute the multiplication table with respect to $\{1,x,y,z\}$. Squares are determined by the trace pairing so it suffices to compute terms such as $xy$. Note
$$tr(x\,xy) = x^2 \tr y = 0 \ , \ \tr(yxy) = 0$$
so $xy$ is in the orthogonal complement of $Kx + Ky$, that is, $xy = cz$ for some scalar $c \in K$. In fact, $\tr xyz = \tr cz^2 = 2cz^2$ and $z^2 \neq 0$ for otherwise, $z$ generates a nilpotent two-sided ideal. Hence, $xy$ and similarly $yz,zx$ are  determined by the trace pairing and $\tr xyz$. 

Given the trace pairing, we observe that there are only two possible values for $\tr xyz$ as follows. 
$$ c^2z^2 = xyxy = -x^2y^2.$$
Hence there are exactly two solutions for $c$ and hence for $\tr xyz$. Moreover, the two solutions are negatives of each other. Since $x,y,z$ skew commute, we see now that $A,A'$ are isomorphic if $\tr xyz$ is the same scalar in both cases while $A,A'$ are anti-isomorphic if $\tr xyz$ are negatives of each other. 
\end{pf}

\begin{corollary}
Assume that $\Pic Z$ has no 2-torsion. Let $A$ be a quaternion $\O_Z$-order with trace pairing $P:sA \rtar (sA)^*$. Then $A \simeq Cl_0(Q)$ where $Q = \sqrt{-1} \Adj P (\det P)^{-\frac{1}{2}}$.
\end{corollary}
\begin{pf}
From proposition~\ref{pfindQ}, $Cl_0(Q)$ has the same trace pairing as $A$. We thus know by the previous lemma that $Cl_0(Q)$ is isomorphic to $A$ or $A^{\mbox{op}}$. However, from the relations of the Clifford algebra, we see immediately that $Cl_0(Q)$ is isomorphic to its opposite algebra so the corollary is proven. 
\end{pf}

\end{section}

%\begin{section}{KILL this section???? } \label{soldbsv}  \marginpar{soldbsv}

%\marginpar{sentence added}
%I think this section should be joined with the next.

%In this section we review some facts about 
%Brauer-Severi varieties of Clifford algebras.  Firstly
%we discuss what happens for a nondegenerate form over
%an algebraically closed field $k$.
%Recall that the projectived spinor bundle of a quadric bundle is the 
%natural embedding space of the maximal isotropic Grassmanian
%or spinor variety.
%If the rank of $V$ is $2n$ then 
%$$V(\Cl_0(Q)) \simeq \P(\wedge^{\even} W) \cup
%\P(\wedge^{\odd} W) \supset OGr(n,2n) \cup OGr(n,2n).$$  It is common
%to take one connected component, which we may not do since they will
%meet when $Q$ degenerates.
%If the rank of $V=2n+1$ is odd, then  
%$V(\Cl_0(Q)) \simeq \P(\wedge^\bullet W) \supset OGr(n,2n-1).$
%This describes the generic fibre of a bundle of quadrics.

%\begin{proposition}
%The Brauer-Severi variety $V(\Cl_0(Q))$ is the projectived spinor bundle
%of the quadric bundle $Q$, if $r=2n+1$ is odd it is a generic $\P^N$
%bundle on $Z$ with 
%$$ N = 2^{n}-1.$$
%If $r=2n$ is even then it is a generic $\P^N$ bundle on a double
%cover of $Z$ with 
%$$ N = 2^{n-1}-1.$$
%\end{proposition}

%At this point we include some good local computations which will
%be useful later.  In particular for 
%$$Q=\begin{pmatrix}
%1 & 0 & 0 \\
%0 & a & b \\
%0 & b & c
%\end{pmatrix}$$
%we know that $V(\Cl_0(Q)) \simeq V(Z^2+aX^2 + b XY +cY^2=0)$
%inside where $X,Y,Z$ are the coordinates of $\P^2$.

%\end{section}

\begin{section}{Brauer-Severi Varieties of Even Clifford Algebras} \label{sbsvofcliff}

In this section we consider a quadratic form $Q: V \otimes V \rtar \L$ where $V$ is a rank 3 vector bundle on a smooth variety $Z$. 
For an algebra of rank $n^2$ over a scheme $Z$ we write $BS(A)$
for the Brauer-Severi scheme ${\rm BSev}_n(A,\O_Z)$ as defined in \cite{VdBS}.  
So we write the Brauer-Severi variety of the even Clifford algebra $Cl_0(Q)$ by $BS(Cl_0(Q))$. Recall that from \S\ref{scbundle} that we may view $Q \in H^0(Z,\Sym^2 V^* \otimes \L)$ and its zero locus $X=X(Q)$ is a conic bundle in $\P(V^*)$. The objective of this section is to show that $X(Q) = BS(Cl_0(Q))$. Hence the maps relating quadratic forms to quaternion algebras and Brauer-Severi varieties are compatible.

The first task is to show that the Brauer-Severi variety naturally embeds in $\P(V^*)$. We can and will assume that $V$ has been normalized as in section~\S\ref{scbundle} so $\L=\det V$. 

\begin{proposition} \label{pacembed}  %\marginpar{pacembed}
Let $A$ be a quaternion algebra on a smooth variety $Z$. Then there is a closed embedding of $BS(A)$ into $\P(A/\O_Z)$. This map sends a codimension two ideal $I < A$ to the one codimensional subsheaf $I + \O_Z < A/\O_Z$. (Here codimension is as of vector spaces over $k$). 
\end{proposition}
\begin{pf}
Recall that $\P(A/\O_Z)$ is the fine moduli space parametrizing one codimensional subsheaves of $A$ which contain $\O_Z$. A map $BS(A) \rtar \P(A/\O_Z)$ can thus be constructed functorially as follows. Let $f:T \rtar Z$ be a test scheme and $I < f^* A$ a left ideal such that $f^*A/I$ is flat over $T$ of constant rank 2. We seek to show that $f^*A/(I + \O_T)$ is a line bundle on $T$ which will give our required map $BS(A) \rtar \P(A/\O_Z)$. To this end, we may assume $\O_T$ is local with maximal ideal $\m$ and we need to show $\Tor^T_1(\O_T/\m,  f^*A/(I + \O_T)) = 0$. Flatness of $f^*A/I$ gives an exact sequence
$$ 0 \rtar \Tor^T_1(\O_T/\m, f^*A/(I + \O_T)) \rtar \O_T/\m \otimes \O_T/\O_T \cap I \rtar \O_T/\m \otimes f^*A/I  .$$
It suffices to show that the map on the right is injective, which, since $\O_T/\m \otimes \O_T/\O_T \cap I \simeq \O_T/\m$ is simple, fails precisely when $\O_T \subset \m f^*A + I$. Suppose this occurs. Now $\m f^*A + I$ is an ideal containing $1$ so must be $f^*A$. Nakayama's lemma now implies that $I = f^*A$, a contradiction. We conclude that $f^*A/(I + \O_T)$ is a line bundle on $T$ so our map $i:BS(A) \rtar \P(A/\O_Z)$ is well defined. 

Now $BS(A)$ is projective over $Z$. We show that our map $i$ is an embedding by showing that it separates points and tangent vectors. This is clear if the points lie over different points of $Z$ or the tangent vector is horizontal. We can thus restrict our attention to some closed fibre $A_0$ of $A$. Let $I_1,I_2$ be distinct two-dimensional ideals in $A_0$. If they are not separated by $i$, then $I_1+k = I_2 + k$. It follows that the ideal $I_1+I_2 = I_1 + k$ which gives a contradiction since the only ideal containing $k$ is $A_0$. 

Now let $k[\e]$ be the ring of dual numbers and $I_1,I_2 < A_0 \otimes k[\e]$ be ideals which are flat over $k[\e]$. 
%\marginpar{tidy proof}
They correspond to vertical tangent vectors in the Brauer-Severi variety which we will assume to be distinct. If they are not separated by $i$ then $I_1 + k[\e] = I_2 + k[\e]$. As in the previous case, $I_1 + I_2 \subset I_1 + k[\e]$ and a contradiction arises unless $I_1 + I_2 = I_1 + \e k[\e]$. Now flatness of $I_1$ implies that $\e A_0 \cap (I_1 + \e k[\e])$ is a 3-dimensional $A_0$-module containing $\e$. However, $A_0 \e$ is already 4 dimensional so we obtain a contradiction once more. 
\end{pf}

\begin{theorem}  \label{tcommute}  %\marginpar{tcommute}
Consider a quadratic form $Q$ on a rank 3 vector bundle $V$ on a smooth variety $Z$ as above. Then $BS(Cl_0(Q)) = X(Q) \subset \P(V^*)$. 
\end{theorem}
\begin{pf}
We carry out the computation at the universal closed point. Hence $V$ is a vector space say with basis $x,y,z$ and $Q$ is given by a $3\times 3$-matrix $(q_{ij})$ with entries in $k$. The even Clifford algebra $A:= Cl_0(Q)$ has basis 
$$Z := x \wedge y , X:= y \wedge z , Y:= z \wedge x, 1 .$$
Recall $A/k \simeq V^*$. We will write elements of $A^*$ as row vectors with respect to the basis dual to $1,X,Y,Z$. This means an element $\a \in (A/k)^*\simeq V$ has the form $\a = (0 \ \a_1 \ \a_2 \ \a_3)$. We compute the closed condition for $\a$ to be in the image of $i: BS(A) \rtar \P(V^*)$. Note that $\ker \a < A$ is 3-dimensional and it is in the image of $i$ precisely when the maximal left ideal $I$ in $\ker \a$ is two dimensional. But using the right $A$-module structure on $A^*$ we can write 
$$I = \ker \a \cap \ker (\a X) \cap \ker (\a Y) \cap \ker (\a Z) .$$
Now a short computation shows $\ker \a, \ker \a X, \ker \a Y, \ker \a Z$ are the rows of the matrix below $M:=$
$$\begin{pmatrix}
0 & \a_1 & \a_2 & \a_3 \\
\a_1 & 2 q_{23} \a_1 & -q_{33} \a_3 & 2 q_{12} \a_1 + q_{22}\a_2 + 2 q_{23} \a_3 \\
\a_2 & 2q_{13} \a_1 + 2q_{23}\a_2 + q_{33}\a_3 & 2q_{13}\a_2 & -q_{11}\a_1 \\
\a_3 & -q_{22}\a_2 & q_{11}\a_1 + 2q_{12}\a_2 + 2q_{13}\a_3 & 2q_{12}\a_3
\end{pmatrix}$$
Also $I = \ker M$ is two-dimensional precisely when $M$ has rank two, that is, all $3\times 3$-minors vanish. Let 
$$q(\a_1,\a_2,\a_3) = q_{11}\a_1^2 + q_{22}\a_2^2 + q_{33}\a_3^2 + 2 q_{12}\a_1\a_2 + 2q_{13}\a_1\a_3 + 2 q_{23}\a_2\a_3.$$  
Then all $3 \times 3$-minors are multiples of $q$ and furthermore, the $(4,3),(3,2),(2,4)$ minors are $\a_1 q, -\a_2 q, \a_3 q$ respectively. Hence the closed condition for $\a$ to be in the Brauer-Severi variety is $q(\a) = 0$. This proves that indeed $BS(Cl_0(Q)) = X(Q)$. 
\end{pf}

\end{section}

\begin{section}{Quaternion algebras of conic bundles}  \label{salgconic}   %\marginpar{salgconic}

In this section, we give a direct method for recovering quaternion algebras from their Brauer-Severi variety. Let $\pi:X \rightarrow Z$ be a conic bundle on a smooth variety $Z$.  So $X$ is a 
Gorenstein variety and the relative dualizing sheaf $\omega_{X/Z} = \omega_X \otimes \pi^* \omega_Z^{-1}$ is a line bundle.  We need the following facts
\begin{lemma}
We have natural isomorphisms
$$ \pi_* \omega_{X/Z} =0$$
$$ R^1\pi_* \omega_{X/Z} = \O_Z$$
$$ H^{i+1}(X,\omega_{X/Z}) = H^i(Z,\O_Z) .$$
\end{lemma}
\begin{pf}
The first statements follow from the fact that $R\pi_*\O_X \simeq \O_Z,$
and the last line follows from the Leray spectral sequence.
\end{pf}
\begin{definition} \label{ldJ}   %\marginpar{dJ}
We define a rank two vector bundle $J$ on $X$ as follows. From the above lemma we see that $H^1(X,\omega_{X/Z}) = H^0(Z,\O_Z)$. Hence $1 \in H^0(Z,\O_Z)$ determines an extension 
$$ 0 \rightarrow \omega_{X/Z} \rightarrow J^* \rightarrow \O_X \rightarrow 0.$$
It is essentially unique. We call this extension the {\it Euler sequence} of the conic bundle and $J$ is the {\it dual Euler extension}.
\end{definition}

The terminology derives from the fact that if we restrict to smooth fibres (or pull back to an \'etale
cover of the locus  of smooth fibres) we obtain
 the usual Euler sequence
$$ 0 \rightarrow \O(-2) \rightarrow \O(-1)^2 \rightarrow \O \rightarrow 0.$$

\begin{lemma}  \label{lrpijstar}  %\marginpar{lrpijstar}
We have the following natural isomorphism
$$ R\pi_* J^* =0.$$
\end{lemma}
\begin{pf}
This follows from the long exact sequence formed
on pushing down the extension above, together with the 
natural isomorphism $R^i\pi_* \O_Z \simeq R^{i+1}\pi_* \omega_{X/Z}.$
\end{pf}

We are primarily interested in the dual of the previous exact
sequence and its pushforward to $Z$.

$$ 0 \rightarrow \O_X \rightarrow J \rightarrow T_{X/Z} \rightarrow 0,$$
where $T_{X/Z}$ is the relative tangent sheaf.
Let $A = \pi_* \cEnd_X(J).$
We have the following result.
\begin{proposition} \label{ppiJ}  %\marginpar{ppiJ}
There is an isomorphism of sheaves on $Z$,
$$ A \simeq \pi_*J.$$
\end{proposition}
\begin{pf}
Apply $-\otimes J$ to the sequence
$$0 \rightarrow \omega_{X/Z} \rightarrow J^* \rightarrow \O_Z \rightarrow 0,$$
to obtain 
$$0 \rightarrow \omega_{X/Z} \otimes J \rightarrow J^*\otimes J \rightarrow J \rightarrow 0.$$
By relative duality we know that $R\pi_*(\omega_{X/Z} \otimes J )$
is dual to $R\pi_* J^* =0$, so the result follows on pushing forward
to $Z$.
\end{pf}

\begin{corollary}  \label{cpushfwdexact}  %\marginpar{cpushfwdexact}
If we push forward the 
exact sequence $$ 0 \rightarrow \O_X \rightarrow J \rightarrow
T_{X/Z} \rightarrow 0$$ we obtain  
$$ 0 \rightarrow \O_Z \rightarrow A \rightarrow A/\O_Z \rightarrow 0$$
and we have a natural isomorphisms
$$ A/\O_Z \simeq \pi_* T_{X/Z} \simeq \pi_* \omega^{-1}_{X/Z} $$ $$
\simeq \pi_* (\O_X(-K_X) \otimes \pi^*\O_Z(K_Z)) \simeq 
\pi_* \O_X(-K_X)  \otimes \O_Z(K_Z).$$
So we see that $A/\O$ is the pushforward of a line bundle.
\end{corollary}

\begin{proposition}  \label{pendisquat}  %\marginpar{pendisquat}
The algebra $A = \pi_* \cEnd_X J$ is quaternion. 
\end{proposition}
\begin{pf}
Lemma~\ref{lrantican} shows that $\pi_* \omega_{X/Z}^{-1}$ is a rank three vector bundle so the exact sequence in the previous corollary reveals that $A$ is locally free of rank four. Now $\cEnd_X J$ is quaternion since it is Azumaya so we may push forward the trace map to obtain a trace map $\mbox{tr}:A \rtar \O_Z$. The conditions on tr for $A$ to be a quaternion algebra are inherited from the corresponding conditions on $\cEnd_X J$. 
\end{pf}

We wish to show that under mild hypotheses, conic bundles and quaternion algebras are in bijective correspondence under the maps 
$$\{\pi:X\rtar Z\}  \mapsto \pi_* \cEnd_X J, \hspace{1cm} A \mapsto BS(A).$$

Under this correspondence, we obtain another important interpretation of $J$. 
Let $A$ be a locally Clifford algebra over $Z$ of rank 4 so that the Brauer-Severi variety $\pi:BS(A) \rightarrow Z$ is a conic bundle by theorem~\ref{tcommute}. Since $BS(A)$ parametrizes two dimensional cyclic representations of $A$ there is a universal cyclic representation $J$ with natural maps $\pi^*A \rightarrow J \rightarrow 0$ and $\pi^*A \rightarrow \cEnd(J)$. We will show that this $J$ corresponds to the one obtained from the conic bundle $BS(A)$ via the Euler sequence. 

We start with a conic bundle $\pi:X \rtar Z$ and seek to show, under some hypotheses, that $X$ is naturally isomorphic to $BS(\pi_* \cEnd_X J)$. The following proposition is the first step.

\begin{proposition}  \label{pmaptoBS}   %\marginpar{pmaptoBS}
Consider the map in the Euler sequence $J^* \rtar \O_Z$ and the induced quotient map $q:\cEnd_X J \rtar J$. 
\begin{enumerate}
\item The composed map 
$$ p:\pi^* \pi_* \cEnd_X J \rtar \cEnd_X J \xrightarrow{q} J $$
is a surjective map of $\pi^* \pi_* \cEnd_X J$-modules. It naturally induces a morphism of varieties $\phi:X \rtar BS(\pi_*\cEnd_X J)$. 
\item The surjection $\pi^*\pi_*\omega_{X/Z}^{-1} \rtar \omega_{X/Z}^{-1}$ defines a map $\psi:X \rtar \P(\pi_*\omega_{X/Z}^{-1})$ and this maps $\psi$ and $\phi$ are compatible with the map $BS(\pi_*\cEnd_X J) \rtar \P(\pi_*\omega_{X/Z}^{-1})$ defined in proposition~\ref{pacembed}. 
\end{enumerate}
\end{proposition}
\begin{pf}
First observe that $q$ is a morphism of $\cEnd_X J$-modules so $p$ is a morphism of $\pi^* \pi_* \cEnd_X J$-modules. To prove 1), it remains only to show that $p$ is surjective since $J$ is flat over $X$ of constant rank two. 
Recall the exact sequence 
$$ 0 \rtar \pi_* \O_X \rtar \pi_* J \rtar \pi_* \omega_{X/Z}^{-1} \rtar 0 .$$
We may pull this back via $\pi$ to obtain a commutative diagram with exact rows
$$
\begin{CD}
L_1\pi^*\pi_* \omega_{X/Z}^{-1}  @>>> \pi^*\pi_* \O_X @>>>   \pi^*\pi_* J @>>>  \pi^*\pi_* \omega_{X/Z}^{-1} @>>> 0  \\
 @VVV @VVV @VVV @VVV \\
0 @>>> \O_X @>>> J @>>> \omega_{X/Z}^{-1} @>>> 0
\end{CD}
$$
Now $\pi^*\pi_* \O_X \rtar \O_X$ is surjective as is $\pi^*\pi_* \omega_{X/Z}^{-1} \rtar \omega_{X/Z}^{-1}$ since $\omega_{X/Z}^{-1}$ is relatively very ample with respect to $\pi$. From proposition~\ref{ppiJ}, we see that $p:\pi^*\pi_*\cEnd_X J = \pi^*\pi_* J \rtar J$ is surjective too and 1) follows. The above commutative diagram also shows that the map $\psi$ is well-defined and compatible with the map in proposition~\ref{pacembed}.
\end{pf}

\begin{theorem}  \label{tBSofend}  %\marginpar{tBSofend}
Let $\pi:X \rtar Z$ be a flat conic bundle and $A = \pi_* \cEnd_X J$. Then the map $\phi:X \rtar BS(A)$ of $Z$-schemes constructed in proposition~\ref{pmaptoBS} is an isomorphism and $J$ is the universal cyclic representation of rank two. Finally, $A$ is locally Clifford of rank 4.
\end{theorem}
\begin{pf}
We know from proposition~\ref{pmaptoBS}, that $\phi$ is compatible with the natural embeddings of $X$ and $BS(A)$ into $\mathbb{P}(\pi_*\omega_{X/Z}^{-1})$. Hence to show it is an isomorphism, it suffices to show that it is an isomorphism on each fibre. Observe that at a closed point $z \in Z$, the Brauer-Severi variety above $z$ is just $BS(A \otimes_Z k(z))$. To compute $A \otimes_Z k(z)$ note that $(\cEnd_X J) \otimes_Z k(z) =  \cEnd_{X_z} (J \otimes_Z k(z))$. Now by construction $J \otimes_Z k(z)$ is the dual Euler extension corresponding to the conic $X_z$ so proposition~\ref{pendisquat} shows that $\End_X (J \otimes_Z k(z))$ is always 4-dimensional. Flatness now gives the base-change condition for $\cEnd_X J$ with respect to $\pi$. Our computation is thus reduced to one on closed fibres.

The isomorphism on closed fibres will follow from the three lemmas below which show the correspondence 
$$\{\pi:X\rtar Z\}  \mapsto \pi_* \cEnd_X J, \hspace{1cm} A \mapsto BS(A) $$
holds on closed fibres. Note that as $\pi$ is flat, there are only three possible fibres, the smooth conic isomorphic to $\P^1$, the pair of lines crossing in a node and finally, the double line. There will be a lemma for each of these cases.

\begin{lemma}  \label{lP1case}  %\marginpar{lP1case}
Let $X = \P^1$ and $J = \O(1) \oplus \O(1)$. Then the dual Euler sequence is 
$$  0 \rtar \O_X \rtar J \rtar \O(2) \rtar 0  $$
and $A = \End_X J$ is the full $2 \times 2$-matrix algebra over $k$. The map $\phi:X \rtar BS(A)$ of proposition~\ref{pmaptoBS} is an isomorphism. Furthermore, the map $p:A \otimes_k \O_X \rtar J$ of that proposition exhibits $J$ as the universal cyclic representation of $A$ of rank two.
\end{lemma}
\begin{pf}
We omit the proof of this easy fact, most of which is well-known.
\end{pf}

\begin{lemma}  \label{lXcase}  %\marginpar{lXcase}
Let $X$ be the union of two distinct lines $l,l'$ in $\P^2$. Let $p,p'$ be points on $l,l'$ respectively which are not nodal.   Then the dual Euler sequence is 
$$  0 \rtar \O_X \rtar \O(p) \oplus \O(p') \rtar \O(p+p') \rtar 0  $$
Setting $J = \O(p) \oplus \O(p')$ we have $A = \End_X J$ is the algebra~2) in theorem~\ref{tclassquat}. The map $\phi:X \rtar BS(A)$ of proposition~\ref{pmaptoBS} is an isomorphism. Furthermore, the map $p:A \otimes_k \O_X \rtar J$ of that proposition exhibits $J$ as the universal cyclic representation of $A$ of rank two.
\end{lemma}
\begin{pf}
The Euler sequence above is clear. Now the Clifford algebra~2) of theorem~\ref{tclassquat} has a Peirce decomposition which allows it to be written schematically as
$$A' :=\begin{pmatrix}
  k & k\e \\
  k\e & k
  \end{pmatrix}, \hspace{1cm} \mbox{where}\ \e^2 = 0.
  $$
Now 
$$\End_X J = 
\begin{pmatrix}
\Hom_X(\O(p),\O(p)) & \Hom_X(\O(p'),\O(p)) \\ 
 \Hom_X(\O(p),\O(p')) & \Hom_X(\O(p'),\O(p')) 
\end{pmatrix}$$
and the algebra isomorphism $A \simeq A'$ is easily obtained by matching up the two Peirce decompositions. It is well-known that $BS(A)$ is isomorphic to $X$ (as can be determined using theorem~\ref{tcommute} for example) from which one easily observes that $\phi$ is in isomorphism and $J$ is the universal cyclic representation.
\end{pf}

\begin{lemma}  \label{l2linecase}  %\marginpar{l2linecase}
Let $R = k[u,v,w]/(w^2)$ and $X \subset \P^2$ be the double line $\Proj R$. Let $A=k\langle x,y\rangle$ be the algebra~3) of theorem~\ref{tclassquat}. Let $M$ be the graded $A \otimes_k R$-module 
$$  M := A \otimes_k R/(R(w+vx-uy) + R(wx-uxy) + R(-wy+vxy) + Rwxy) $$
and $J$ be the corresponding sheaf on $X$. Then the dual Euler sequence is 
$$  0 \rtar \O_X \rtar J \rtar \O_X(1) \rtar 0  $$
and $A \simeq \End_X J$. The map $\phi:X \rtar BS(A)$ of proposition~\ref{pmaptoBS} is an isomorphism. Furthermore, the map $p:A \otimes_k \O_X \rtar J$ of that proposition exhibits $J$ as the universal cyclic representation of $A$ of rank two.
\end{lemma}
\begin{pf}
Since the dual Euler extension is the essentially unique non-split extension of $\O_X(1)$ by $\O_X$, verifying the Euler sequence amounts to showing that the cokernel $N$ of $R \rtar M: r \mapsto 1 \otimes r$ is a Serre module for $\O_X(1)$. Now 
$ R_{>0} wxy \subset R(wx-uxy) + R(-wy+vxy)$ so up to a finite dimensional vector space we have 
$$ N = \frac{Rx \oplus Ry \oplus Rxy}{R(vx-uy) + R(wx-uxy) + R(-wy+vxy)}  .$$
But the Koszul complex for $k[u,v,w]$ shows that this is indeed a Serre module for $\O_X(1)$. 

Since $M$ is an $A$-module we certainly have $A \subset \End_X J$. But proposition~\ref{pendisquat} shows that $\End_X J$ is 4-dimensional so we have equality. We know $BS(A)$ is the double line so it follows that $\phi$ must be an isomorphism and $J$ is the universal cyclic representation.
\end{pf}

Proposition~\ref{pendisquat} shows that $A$ is quaternion while the fibre-wise computations above show that all the closed fibres of $A$ are generated as a $k$-algebra by two elements. Proposition~\ref{ploccliff} now ensures that $A$ is locally Clifford of rank 4. This completes the proof of the theorem. 
%\marginpar{tighten proof}
\end{pf}

\begin{theorem}  \label{tendofBS} %\marginpar{tendofBS}
Let $A$ be a locally Clifford algebra over $Z$ of rank four and $\pi:X=BS(A) \rtar Z$ be the Brauer-Severi variety. Then $\pi$ is a flat conic bundle and $A \simeq \pi_* \cEnd_X J$ where $J$ is universal cyclic representation of rank two. Furthermore, $J$ is the dual Euler extension associated to the conic bundle $\pi:X \rtar Z$. Consequently, there is a bijection between flat conic bundles and locally Clifford algebras of rank four.
\end{theorem}
\begin{pf}
Since $A$ is locally Clifford, it is locally even Clifford so $BS(A)$ is a conic bundle by theorem~\ref{tcommute}. None of the fibres of $BS(A)$ are $\P^2$ so it is in fact a flat conic bundle. We have by definition of universal representation a surjective module map $\pi^* A \rtar J$ and a map $\pi^*A \rtar \cEnd J$. Hence there is an algebra map $A \rtar \pi_* \cEnd_X J$. It is an isomorphism by the fibre-wise computations in lemmas~\ref{lP1case},\ref{lXcase} and \ref{l2linecase}. 

The fibre-wise computation also shows that on every closed fibre $X_z$ for $z \in Z$, we have a non-split sequence 
$$ 0 \rtar \O_{X_z} \rtar J|_{X_z}  \rtar \omega_{X_z}^{-1}  \rtar 0  .$$
This shows that $T := J/\O_X \simeq \omega_{X/Z}^{-1} \otimes_Z \pi^* \M$ for some line bundle $\M \in \Pic Z$. We need to show that $\M \simeq \O_Z$. Now $R^1\pi_* T^* = R^1\pi_* \omega_{X/Z} \otimes \M^* = \M^*$ so it suffices to show that $R^1\pi_* T^* \simeq \O_Z$. 

Note that $R\Gamma (X_z,J^*|_{X_z}) = 0$ by lemma~\ref{lrpijstar} so $R\pi_* J^* = 0$ too. Consider the universal ideal $I$ and the exact sequence 
$$  0 \rtar I \rtar \pi^* A \rtar J \rtar 0 .$$
We dualize to obtain a commutative diagram with exact rows.
$$
\begin{CD}
 0 @>>> T^* @>>> \pi^*(A/\O_Z)^* @>>> I^* @>>> 0 \\
 @VVV  @VVV @VVV @|  @VVV  \\
0 @>>> J^* @>>> \pi^*A^* @>{\a}>> I^* @>>> 0 
\end{CD}.
$$
Now $R\pi_* J^* = 0$ so $\pi_* \a: A^* \rtar \pi_*I^*$ is an isomorphism. Hence, we see that 
$$R^1\pi_* T^* = \mbox{coker} ( (A/\O_Z)^* \rtar \pi_* I^*) = \mbox{coker} ( (A/\O_Z)^* \rtar A^*) = \O_Z.$$
This completes the proof.
\end{pf}

\end{section}

\begin{section}{Chern classes and {\bf $-K^3$}}  \label{sc2}

We will now compare some other invariants of
conic bundles and quaternion orders.  In this section we assume
that $X$ is a smooth threefold which is a conic bundle over a smooth surface
$Z$.
Riemann-Roch gives us the following formula
for any coherent sheaf $\E$ on a smooth threefold $X$,
$$\chi(\E)=\deg(\ch(\E).\td(T_X))_3.$$
We will temporarily write $c_i$ as shorthand for the
Chern classes of the tangent bundle $c_i(T_X)$.
We will write 
$$c_3 = \chi_{\top}(X).$$
simply as notation.
Using Riemann-Roch gives
$$\chi(\O_X) = \frac{1}{24} c_1c_2.$$
Also 
$$c_1=-K_X.$$

Now applying Riemann-Roch again gives 
\begin{eqnarray*}
\chi(T_X) & = & \deg(\ch(T_X).\td(T_X))_3 \\
 & = & \frac{1}{6}(c_1^3-3c_1c_2+3c_3)+\frac{1}{12}(c_1^2+c_2)c_1
+\frac{1}{4}c_1(c_1^2-2c_2)+\frac{1}{8}c_1c_2 \\
&  = & \frac{1}{2} c_1^3- \frac{19}{24} c_1c_2 +\frac{1}{2} c_3\\
& = & -\frac{1}{2} K_X^3 -19 \chi(\O_X) + \frac{1}{2} \chi_{\top} (X).
\end{eqnarray*}
So we now have
$$ -\frac{K_X^3}{2} + \frac{\chi_{\top}(X)}{2} = \chi(T_X) + 19 \chi(\O_X).$$

For a standard conic bundle $\pi:X \rightarrow Z$
we can simplify the formulas further.  Let the discriminant be $D$, 
we have that 
$$\chi_{\top}(X) = 2 \chi_{\top}(Z) + \chi_{\top}(D),$$
$$\chi(\O_X) = \chi(\O_Z)$$
as in \cite{IP} Lemma 7.1.10.
For a standard conic bundle $\pi:X \rightarrow Z$ with discriminant $D,$
 we also have the following exact sequence which 
follows from local computations
$$ 0 \rightarrow T_{X/Z} \rightarrow T_X \rightarrow \pi^* T_Z
\rightarrow i_*N_{Z/D} \rightarrow 0,$$
where $N_{Z/D}$ is the normal bundle of $D$ in $Z$, and $i$ is the isomorphism
from $D$ to the singular locus of $\pi^{-1}(D)$.
So we can compute 
$$-K_X^3/2 = 19\chi(\O_Z) + \chi(T_Z) 
- \chi_{\top}(Z)-\chi_{\top}(D)/2 -\chi(\O_D(D)) + \chi(A/\O).$$
$$= \frac{11}{4} K^2 -\frac{1}{4} \chi_{\top}(Z)  +K_Z.D+ \chi(A/\O)$$
$$ = 3K_Z^2-3 \chi(\O_Z) +K_Z.D +\chi(A/\O).$$
So we obtain the following result.
\begin{proposition}
Let $\pi:X \rightarrow Z$ be a standard 
conic bundle with associated quaternion order
$A$ and discriminant $D$.
Then
$$-K_X^3=  6K_Z^2+3K_Z.D+D^2 -c_2(A).$$
% $$H^i(Z,\pi_*T_{X/Z}) = H^i(Z,A/\O_Z)$$
\end{proposition}
\begin{pf}
We use the fact that $c_1(A)=-D$ and Riemann-Roch for surfaces
with the above computation.,
\end{pf}

If we restrict to the case where $Z = \P^2$ 
and let $\deg D =d$ then we get the formulas 
$$-K_X^3=48-6d+2\chi(A/O)$$
$$ -K^3_X = 54-9d+d^2-c_2(A).$$

\end{section}

\begin{section}{Del Pezzo Orders and Conic Bundles} \label{sdPO}

\begin{subsection}{Del Pezzo Orders} \label{squadofdp}  %\marginpar{squadofdp}

We are interested studying del Pezzo quaternion orders and their associated conic bundles. The minimal del Pezzo orders were classified in terms of their ramification data $(\widetilde{D} \rtar D \rtar Z)$ in \cite{CK,CI,AdJ}. 
We will only be concerned with minimal terminal quaternion del Pezzo orders.
We will refer to these simply as del Pezzo orders but it should be noted
that there are many other types of del Pezzo orders which are not necessarily
minimal, terminal or quaternion.  Briefly, in the quaternion case, the centre of the order is always $Z=\mathbb{P}^2$, the ramification locus $D\subset Z$ is a nodal curve of degree $d=3,4$ or $5$ and $\widetilde{D}$ is a double cover of $D,$ ramified at the nodes. We denote them by $F^2_d$. For each ramification data, we wish to explicitly construct quadratic forms $Q: \Sym^2 V \rtar \L$ such that the corresponding Clifford algebra $Cl_0(Q)$ has ramification data $(\widetilde{D} \rtar D \rtar Z)$. The centre $Z$ of the del Pezzo order is $\P^2$, so we may use Catanese theory \cite{Cat} to construct $Q$, as has been done by Brown-Corti-Zucconi [BCZ]   We will review that construction.
%(see Appendix \S\ref{asymmres}). 

\begin{proposition}
Let $(\widetilde{D} \rtar D \rtar Z)$ be the ramification data of a
minimal del Pezzo order.  Then the symmetric resolution of 
$L:=\O_{\widetilde{D}}/\O_D$ is one of the following types.
\begin{eqnarray*}
F^2_3: & 0 \rightarrow \O(-2)^3 \rightarrow \O(-1)^3 \rightarrow L \rightarrow 0 \\
F^2_4: & 0 \rightarrow \O(-3)^2 \rightarrow \O(-1)^2 \rightarrow L \rightarrow 0 \\
F^{2}_{5+}: & 0 \rightarrow \O(-3)^5 \rightarrow \O(-2)^5 \rightarrow L \rightarrow 0 \\
F^{2}_{5-}: &  0 \rightarrow \O(-4)\oplus\O(-3)^2 \rightarrow \O(-2)^2 \oplus \O(-1) \rightarrow L \rightarrow 0 \\
\end{eqnarray*}
\end{proposition}
\begin{pf}
Write $\O_{\widetilde{D}} = \O_D \oplus L$ for some 2-torsion line bundle $L$ on $D$.We can resolve the module $\Gamma(L) = \oplus H^0(\P^2,L(i))$ over the
homogeneous coordinate ring of $\P^2$.  This will give a resolution by 
sums of line bundles.  So we may use 
apply Catanese theory \cite{Cat}, which 
requires locally free resolutions of $L$.
The types above all follow from Riemann-Roch calculations.  We 
will work out the case $F^2_4$ in detail and explain why there are two
separate cases for $F^2_5$.

If $\deg D=4$, and $L^{\otimes 2}\simeq \O_D$ is non-trivial. Then 
$h^0(L)=0$ and $\deg L=0$.  We see that $\chi(L(i))=4i-2$ and 
$h^1(L(i)))=h^0(L(-i+K_D))=0$ for $i\geq 2$.  Also
$h^1(L(1))=h^0(L^*)=h^0(L)=0$.  So our resolution begins with $\O(-1)^2$
since $h^0(L(1))=2.$
To find  the required syzygy we twist by one and 
compute $h^0(L(2))=6$ and $h^0(\O(1)^2)=6,$ so no syzygy is required in
this degree.  Twisting once more yields $h^0(L(3))=10$ and $h^0(\O(2)^2)=12,$ 
so we require $\O(-3)^2$ as a syzygy.  Checking Hilbert series 
show that the resolution is complete at this point.

In the case where $D$ is a smooth
 quintic case, $L(1)$ is a
theta characteristic, and its parity effects the Riemann-Roch calculation.  We know
by Clifford's Theorem that $h^0(L(1)) \leq 3$.  So since $L(1)$ has degree
5, if $h^0(\P^2,L(1)) =2$ or 3, then we know that $L(1) \simeq \O_D(1)$
or $L(1) \simeq \O_D(1) \otimes \O_D(p-q)$ by exercise B-1, p.264 of \cite{ACGH}.
The first case is certainly not possible, and in the last case we would
require that $\O(2p)\simeq \O(2q)$ giving that $D$ is hyperelliptic.
By exercise B-2, p.221 loc. cit, we see that this is also impossible.
Hence we have either $h^0(L(1))=0$ or 1, and Riemann-Roch calculations give
the above two resolutions of types $F^{2+}_5$ and $F^{2-}_5$. 
\end{pf}

Since these resolutions are symmetric when obtain quadratic forms.
Two of the resolutions types yield quadratic forms with vector bundles
that do not have rank three.
We make a simple adjustment to the case of $F^2_4$ by adding
an $\O(-2)$ to each rank two vector bundle to obtain the
new resolution:
\begin{eqnarray*}
F^2_4: & 0 \rightarrow \O(-2)\oplus \O(-3)^2 \rightarrow \O(-2) \oplus 
\O(-1)^2 \rightarrow L \rightarrow 0 \\
\end{eqnarray*}  We also need to make a more
complicated adjustment in the quintic even theta characteristic case
$F^{2+}_5$ as explained later.

%So we get 
%\begin{subsection}{$GT^2_7$}
%\marginpar{FIX}
%There is a curve of genus 129 parameterizing the lines in $X_8$,
%the complete intersection of 3 quadrics in $\P^6$.  So I'd
%like to show that the tanget space has dimension one but 
%it doesn't.  Perhaps its not reduced?
%\end{subsection}

The first term of the resolution is the vector bundle $V^*$
and the resolution can be chosen to be symmetric so we have
the map $\widetilde{Q} : V \otimes \L^* \rightarrow V^*$.  So we may
construct the even Clifford algebra.
We can choose numbers $a_1,a_2,a_3,d$ so that $V= \oplus \O(2a_i+d)$.
Then $Q$ can be presented as a symmetric matrix with entries which
are forms of degree $\deg Q_{ij}=a_i+a_j+d$.  The numbers are chosen
to be 

\begin{center}
\begin{tabular}{c|c|c|c|c}
type & $a_1$ & $a_2$ & $a_3$ & $d$ \\
$F^2_3$ & 0  & 0 & 0 & 1 \\
$F^2_4 $ & 0 & 1 & 1 & 0 \\
$F^2_{5-}$ & 0 & 0 & 1 & 1 
\end{tabular}
\end{center}

In this case we can form a homogeneous coordinate ring for 
the even Clifford algebra which is fairly simple.
$$\Gamma(\Cl_0(Q)) = \bigoplus_{i \geq 0} 
H^0(\P^2,\Cl_0(Q) \otimes \O(i))$$

We will form a Clifford algebra $\Cl(Q)$ over the polynomial ring $k[u,v,w]$,
generated by $x_1,x_2,x_3$ with the relations $(x_i,x_j)=Q_{ij}$.
We set the degrees of $x_i$ to be $2a_i+d$ and the degrees
of $u,v,w$ to be 2.  So we get a graded algebra $B$
with the 6 generators $x_1,x_2,x_3,u,v,w$.  
\begin{proposition}  The algebra $\Gamma(\Cl_0(Q))$
is the subalgebra of $\Cl(Q)$ generated by  $u,v,w,x_1x_2,x_2x_3,x_1x_3$.
\end{proposition}
\begin{pf}
We first check that the Hilbert series are the same.
For the algebra $A$ we have that 
$$ A_n = 
H^0(\P^2,(\O \oplus \O(a_1+a_2+d) \oplus \O(a_2+a_3+d) \oplus \O(a_1+a_2+d))
\otimes \O(d+n)).$$
We also see that $B$ is generated as a module over $k[u,v,w]$ 
by $1,x_1x_2,x_2x_3,x_3,x_1$ whose degrees are $0,2a_1+2a_2+2d, 2a_2+2a_3+2d,
2a_3+2a_1+2d$.  So the Hilbert series match.  
Now the construction of the even Clifford algebra shows that
there is a map $\Gamma(\Cl_0(Q)) \rightarrow \Cl(Q)$.
\end{pf}

A similar analysis can be done for the even Clifford algebras
of type $F^{2+}_5$ where we will obtain an algebra $\Cl_0(Q)$ with rank $4^2$.

 Given an order $A$ we define
the Kodaira dimension of $A$ to be give by the growth of
the Hilbert Series of the canonical algebra
$$ \bigoplus H^0(Z,\omega_A^n).$$ 
It is not hard to see that the Kodaira dimension of $A$
is the same as the Kodaira dimension of the associated log
surface as in \cite{CI}.

\begin{proposition}
Let $A$ be an order over smooth surface
$Z$ with $A\otimes k(Z)$ a division algebra,
$H^1(Z,A)=0$, and $\kod(A)=-\infty$.  Then if $B$ is Morita equivalent
to $A$ and has the same rank then $c_2(B) \geq c_2(A)$.
\end{proposition}
\begin{pf}
Since $B$ and $A$ have the same first Chern class since the
have the same rank and discriminant.  So 
Riemann-Roch yields
$\chi(A)-\chi(B)=c_2(B)-c_2(A)$.  Now $h^2(Z,B)=h^0(Z,\omega_B)=0$
since $A$ has $\kod(A)=-\infty$ which is a Morita invariant.
Also $H^0(Z,A) =H^0(Z,B)=k$ since $A$ is in a division algebra.
\end{pf}

In the cases under consideration, if we let $A=\Cl_0(Q)$ then
 $h^i(A/\O_{\P^2}) = h^i(V^*)=0$ we see that there are
no deformations
of $A$ as an order over $Z$, or in other words $A$ is rigid.
The above Proposition also shows that if $A$ 
is of type $F^2_3, F^2_4, F^{2+}_5$ then $A$ has a minimal second
Chern class
among Morita equivalent orders with the same rank.  Results of \cite{AdJ}
show that the moduli space of such orders is a proper scheme
of dimension zero.  We conjecture further that the moduli
space is a single point.

\begin{conjecture}
The even Clifford algebras $\Cl_0(Q)$ are the only orders 
which have the same rank, second Chern class and are 
Morita equivalent to $\Cl_0(Q)$.
\end{conjecture}

We suspect this conjecture is true for all quaternion minimal terminal
del Pezzo orders, but we have less evidence for type $F^{2+}_5$ since
we do not know if the second Chern class is minimal.

%
%\begin{subsection}{$F^2_{2+}$}
%In this case the curve $D$ is quintic and so has genus 5.
%Let $\L$ be a line bundle on $D$ such that $\L^2 \simeq \O_D$,
%and we also supose that the theta characteristic of $\L$
%is even.  This amounts to saying that the degree zero line
%bundle $\L(1)$ has no sections.  So Rieman-Roch calculations
%show that there is a resolution 
%$$ 0 \rightarrow \O(-3)^5 \rightarrow \O(-2)^5 \rightarrow \L \rightarrow 0.$$
%But this is larger than we want, and we see that we may split
%this resolution and write
%$$ 0 \rightarrow \Omega^1(-2) \oplus \O(-3) \rightarrow \Omega^1 \oplus \O(-2)
%\rightarrow \L.$$
%\end{subsection}

\end{subsection}

\begin{subsection}{Conic Bundles of del Pezzo Orders}
We now describe the associated conic bundle of the del Pezzo
orders.  Since each type is significantly different we will
discuss each separately.  We first note that
the conic bundles are all Fano by the following result.

\begin{proposition}
Let $\pi:X \rightarrow Z$ be a standard conic bundle
and suppose that $X=V(Q=0) \subseteq \P(V^*)$ where
$V$ is normalized.  If for any curve $C$ in $Z$
and a surjection $V^* \rightarrow L$ where $L$ is  a
line bundle supported on $C$, we have that $\deg L - K_Z.C >0$
then $X$ is Fano.
\end{proposition}
\begin{pf}
Consider the divisor $-K_{\P(V^*)}-X$ in $\P(V^*)$.  We will show
under the given conditions that $-K_{\P(V^*)}-X$ is an ample divisor on $\P(V^*)$
so the restriction $-K_X=(-K_{\P(V^*)}-X)_{|X}$ will be ample on $X$.
The discussion above Proposition~\ref{gorsch} shows that 
$-K_{X/Z} = H$ and so $-K_{\P(V^*)}-X=H-\pi^*K_Z$ where $H$ is the divisor
on $\P(V^*)$ with the property that $\pi_*\O(H)=V^*$.  
The cone of effective curves of $\P(V^*)$ is generated by fibres and sections
over curves in $Z$.  If $F$ is a fibre then $(H-\pi^*K_Z).F=1$ since $H$
is a section.  If $C$ is a curve and $\sigma:X \rightarrow Z$
is a section, the section $\sigma$ is determined by a quotient line 
bundle $V^*_{|C} \rightarrow L$ on $C$.  One can compute that
$(H-\pi^*K_Z).\sigma(C)=\deg L - K_Z.C$.
\end{pf}

\begin{corollary}
Let $\pi : X \rightarrow \P^2$ be a standard conic bundle in 
$\P(V^*)$ where $V^*$ is normalized.  Then if $\deg M < 3$ for
all lines $L \subset \P^2$ and sub-linebundles $M$ of $V_{|L}$ then
$X$ is Fano.
\end{corollary}
It is a simple matter to verify the above condition for the types 
of interest $F^2_3,F^2_4,F^{2+}_5,F^{2-}_5$.

Since smooth Fano threefold are well described in \cite{IP} and
in particular all the Fano conic bundles are listed.  So we give a
geometric construction of the Brauer-Severi variety of each type of 
del Pezzo orders.

\vspace{2mm} \noindent {$\bf{F^2_3}$:}  In this case the quadratic form $Q$
has rank 3 with linear entries.  So $Q$ can be interpreted as a
polynomial of bidegree $(1,2)$ on $\P^2 \times \P^2$.  The divisor
where the polynomial vanishes 
is the Brauer-Severi variety $X$ with one projection giving a
conic bundle ramified on a cubic and the other presenting $X$
as ruled over $\P^2$.  
So $X=\P(\E)$ where $\E$ is a rank
two vector bundle which is the cokernel of the map $\O(-2) \rightarrow \O^3$
derived from $Q$.

\vspace{2mm} \noindent ${\bf F^2_4}$:  In this case we can consider the quadratic form 
$Q$ to be a polynomial of bidegree $(2,2)$ on $\P^1 \times \P^2$.
The Brauer-Severi variety $X$ is the the double cover of $\P^1 \times \P^2$
ramified on this divisor.  Note that the projection to $\P^1$ presents
$X$ as a quadric surface bundle with 6 singular fibres.

\vspace{2mm} \noindent  ${\bf F^{2+}_5}$: 
%\marginpar{REF appendix \S\ref{sorthoproj}}
These conic bundles may be constructed from a net of quadrics in 
$\P^4$.  The base locus will be a canonically embedded complete 
intersection curve
$C$ of degree 8 and genus 5.  If we choose a point $p \in C$
we can do an orthogonal projection of our quadric bundle.
For each point $q$ in the base $\P^2$ we
by replace the quadric $Q_q$ with the conic $T_p Q_q\cap Q_q$.
Let $\pi(C)$ be the image of the curve $C$ under the projection
$\pi:\P^4-p \rightarrow \P^4$
from $p$.  The curve $\pi(C)\simeq C$ has degree 7 in $\P^3$.
The conic bundle we get from the orthogonal projection is
$X=\Bl_{\pi(C)}\P^3$.  Now we will discuss the the order
associated to this conic bundle.

Let $Q$ be a net of quadrics in $\P^4$ with base
$\P^2$ having coordinates $x,y,z$.  So $Q$ is the vanishing 
locus of $v^tAv$ where $A$ is a symmetric $5 \times 5$ matrix
with entries that are linear in $x,y,z$.  We choose a point
$p$ in the base locus of $Q$ such that $p$ is a smooth 
point of every quadric in the net.  To obtain a conic bundle
we do the standard trick of forming the quadrics $T_p Q_q\cap Q_q$
and taking the image under the projection $\P^4-p \rightarrow \P^3$.
This projection changes $V$ from $\O^5$ to $\O^4$.
The image is a conic bundle which is degenerate on the quintic $\det A=0$.  Let us assume that $p=[0,0,0,0,1]$ and so $A=(a_{ij})$
has $a_{55}=0$.  We can compute that the tangent space
at each point is the vanishing locus $T_p = \{ v \in \P^4: p^T A v \}$,
and $p^TA =(a_{15},a_{25},a_{35},a_{45},0)$.  Since the quadrics
are all smooth at $p$, if we consider $p^TA$ as a $4\times 3$ matrix,
it will have rank three.  So by adjusting bases for both $\P^2$ and $\P^4$
we may assume that $p^TA=(x,y,z,0,0)$.  So this allows us to compute
the structure of the vector bundle $\P(V)$ in $\P^3 \times \P^2$
 above $\P^2$ which contains the conic bundle.  It is given by the
exact sequence 
$$ 0 \rightarrow \O_{\P^2}(-3) \stackrel{p^TA}{\rightarrow} \O_{\P^2}(-2)^4
\rightarrow V \rightarrow 0.$$
So we see that $V \simeq \Omega^1 \oplus \O(-2)$.

So we obtain the new resolution
\begin{eqnarray}
F^{2+}_5: & 0 \rightarrow \Omega^1(-2) \oplus \O(-3) \rightarrow \Omega^1 \oplus \O(-2)
\rightarrow L \rightarrow 0.
\end{eqnarray}

We conjecture that the moduli space of del Pezzo Orders of type $F^{2+}_5$
described above is the curve $C$.
\begin{conjecture}
Given a fixed ramification data of type $F^{2+}_5$, the moduli space
of quaternion orders with fixed Morita equivalence class and Chern classes
is the curve $C$ as constructed above.
\end{conjecture}
There is a corresponding conjecture for conic bundles.
\begin{conjecture}
Let $X=\Bl_{\pi(C)}\P^3$. The moduli space of conic bundles over $\P^2$
which are birational to $X$ over $\P^2$ with fixed anti-canonical degree 
is given by the curve $C$.
\end{conjecture}

\vspace{2mm} \noindent ${\bf F^{2-}_{5}}$:  In this case the Brauer-Severi variety
is the blow up of a cubic threefold along a line $X=\Bl_l V_3$.
To show the relation with $Q$ choose coordinates $u,v,x,y,z$
on $\P^4$ so that the line $l=V(x=y=z=0)$.  Let our cubic
threefold be $V_3=V(f=0)$ and note that since $l \subset V_3$
we have that $f \in (x,y,z)$.  Now and write $f$ as a polynomial in $u,v$
$$f= q_{11} u^2+2q_{12}uv +q_{22}v^2 +2q_{13}u+2q_{23}v+q_{33}.$$
Our quadratic form has entries $Q=(q_{ij})$.
The conic bundle structure is given by the projection from the line.

We can present the geometric version of the conjecture of the
uniqueness of moduli.

\begin{conjecture}
Let $X$ be the Brauer-Severi variety of an order of type 
$F^2_3, F^2_4, F^{2-}_5$ as described in the table below.  Then if $Y$ is 
birational to $X$ over $\P^2$ and has the same anticanonical 
degree then $Y \simeq X$.
\end{conjecture}

For convenience we record some of the results in this section in the
following table.

%\begin{section}{Rational Threefold Conic Bundles}  \label{sratbun}   \marginpar{sratbun}

%The following threefolds are the Brauer-Severi
%varieties of the orders listed below.  They are all
%Fano except $GT^2_7$.  The table
%rows are taken from \cite{IP}.  In the table,
%V_i$ or $V_{i,j}$ indicate divisors of the appropropriate
%multi-degree and a 2 above an arrow indicates that the map is 2:1.

%$CY^2_6$ &  $V_{2,2} \subset \P^2 \times \P^2$ & 12 & 9 \\
%$CY^2_6$ &  $V \stackrel{2}{\rightarrow} \P(\Omega_{\P^2})$ 
%ramified on $-K$ & 12  & 9 \\
%$GT^2_7$ & $\Bl_{\mbox{line}}V_8$ with $V_8 \subset \P^6$ & 2  & 14 \\ 
%$GT^2_8$ & $V \stackrel{2}{\rightarrow} \P^1 \times \P^2$ 
%ramified on $V_{2,4}$ & 6 & 20 \\ 
%\end{tabular}

%Also the variety for $F^2_4$ is also a quadric surface bundle over
%$\P^1$ which degenerates to a cone at six points in $\P^1.$
%This gives us a hyperelliptic curve of genus 2 as the Clifford
%algebra of this quadric bundle, which is a component of the derived
%category.

%Here we list some useful formulas for invariants of standard
%conic bundles.

\begin{tabular}{c|c|c|c|c}
 type & $V^*$ & $BS(A)$ &  $-K_X^3$ & $h^{1,2}$ \\
\hline
$F^2_3$ &  $\O(-1)^3$ & $X_{1,2} \subset \P^2 \times \P^2$ & 30 & 0 \\
$F^2_4$ &  $\O(-2) \oplus \O(-1)^2$ & $X \stackrel{2}{\rightarrow} \P^1 \times \P^2$ 
ramified on $V_{2,2}$ & 24 & 2 \\ 
$F^{2+}_5$ & $\Omega^1\oplus\O(-2)$ & $\Bl_C\P^3, \deg C =7, \g(C) =5$ & 16 & 5 \\
$F^{2-}_5$ &  $\O(-2)^2 \oplus \O(-1)$ & $\Bl_{\mbox{line}}V_3$ with $V_3 \subset \P^3$, & 18 & 5 
\end{tabular}

\end{subsection}
\end{section}

\end{document}